\documentclass[reqno]{amsart}
\usepackage[top=25truemm,bottom=25truemm,left=25truemm,right=25truemm]{geometry}
\usepackage[dvipdfmx]{graphicx}
\usepackage{amsmath,mathrsfs,amssymb}
\usepackage{enumerate}
\usepackage{here}
\usepackage{cases}
\usepackage{color}
\usepackage{bm}
\usepackage{amsmath,amssymb}
\usepackage{mathrsfs}
\usepackage{enumerate}
\usepackage{multirow}
\usepackage{cases}
\usepackage{color}

\newcommand{\editcolor}[1]{
	\if\editnum0#1\fi%
	\if\editnum1\textcolor{white}{#1}\fi}

\newcommand{\rmn}[1]{{\rm #1}}

\newcommand{\ds}{\displaystyle}
\newcommand{\nn}{\nonumber}
\newcommand{\ra}{\rightarrow}

\newcommand{\fa}{\forall}
\newcommand{\ex}{\exists}

\newcommand{\st}{\mbox{~s.t.~}}
\newcommand{\leftright}[3]{\left#1#2\right#3}

\newcommand{\BrS}[1]{\leftright{(}{#1}{)}}
\newcommand{\BrM}[1]{\left\{#1\right\}}

\newcommand{\BrA}[1]{\left|#1\right|}
\newcommand{\BrN}[1]{\left\|#1\right\|}

\newcommand{\BrANd}[1]{|#1|}
\newcommand{\BrNNd}[1]{\|#1\|}
\newcommand{\Ol}[1]{\overline{#1}}

\newcommand{\citerm}[1]{{\rm \cite{#1}}}

\newcommand{\diam}{\textrm{diam}}

\newcommand{\gausssymbol}[1]{\lfloor #1 \rfloor}
\newcommand{\Dim}{d}
\newcommand{\midd}{\;\middle|\;}

\newcommand{\Set}[2]{\left\{#1;\,#2\right\}}
\newcommand{\SetNd}[2]{\{#1;\,#2\}}

\newcommand{\Meas}[1]{\BrA{#1}} 
\newcommand{\supp}{\mbox{supp}} 
\newcommand{\Norm}[2]{\left\|#1\right\|_{#2}} 
\newcommand{\OBallCher}{B}
\newcommand{\OBall}[2]{\OBallCher(#1,#2)}

\newcommand{\Intp}{\Pi}
\newcommand{\Grad}{\nabla}
\newcommand{\Lap}{\Delta}

\newcommand{\Dm}{\Omega}
\newcommand{\DmOl}{\Ol{\Dm}}
\renewcommand{\H}{H}
\newcommand{\DmH}{\Dm_\H}
\newcommand{\DmHOl}{\Ol\Dm_{\H}}

\newcommand{\N}{N}
\newcommand{\PtSetCher}{\mathcal{X}}
\newcommand{\PtSet}{\PtSetCher_{\N}}

\newcommand{\IndexSetCher}{\Lambda}
\newcommand{\IndexSet}[1]{\IndexSetCher(#1)}

\renewcommand{\i}{i}
\renewcommand{\j}{j}
\renewcommand{\k}{k}

\newcommand{\PtChar}{x}
\newcommand{\Pt}[1]{\PtChar_{#1}}

\newcommand{\Dx}{\Delta\x}

\newcommand{\h}{h}
\newcommand{\hN}{\h_{\N}}
\newcommand{\WeightCher}{w}
\newcommand{\Weight}{\WeightCher}

\newcommand{\WeightH}{\Weight_{\h}}

\newcommand{\w}{\Weight}

\newcommand{\wh}{\WeightH}

\newcommand{\wn}[1]{\w^{(#1)}}
\newcommand{\wnh}[1]{\w^{(#1)}_\h}

\newcommand{\FuncCher}{v}
\newcommand{\f}{\FuncCher}

\newcommand{\PvChar}{V}
\newcommand{\Pv}[1]{\PvChar_{#1}}

\newcommand{\PvSetCher}{\mathcal{V}}
\newcommand{\PvSet}{\PvSetCher_{\N}}

\newcommand{\VoroCellChar}{\sigma}
\newcommand{\VoroCell}[1]{\VoroCellChar_{#1}}

\newcommand{\IntpApp}{\Pi_\h}
\newcommand{\GradApp}{\Grad_\h}

\newcommand{\LapApp}{\Delta_\h}

\newcommand{\IndPtSet}{r_{\N}}
\newcommand{\IndPvSet}{d_{\N}}
\newcommand{\IndPvSetLocal}[1]{d_{#1}}

\newcommand{\FsWeightFunc}{\mathcal{W}}
\newcommand{\gconst}{c}

\newcommand{\RegOrder}{m}
\newcommand{\maxi}[1]{\max_{#1=1,2,\dots,\N}}
\newcommand{\WeightOrder}{n}
\newcommand{\WeightOrderSmooth}{k}

\newcommand{\regularconst}{c_0}
\newcommand{\PvTempCher}{\widetilde{\PvChar}}
\newcommand{\PvTemp}[1]{\PvTempCher_{\i}}

\newcommand{\PvDmChar}{\xi}
\newcommand{\PvDmSet}{\Xi}
\newcommand{\PvDm}[1]{\PvDmChar_{#1}}
\newcommand{\PvDmOl}[1]{\Ol{\PvDmChar}_{#1}}

\newcommand{\FsCChar}{C}
\newcommand{\FsCz}[1]{\FsCChar(#1)}
\newcommand{\FsC}[2]{\FsCChar^{#1}(#2)}
\newcommand{\NormCz}[2]{\BrN{#1}_{\FsCz{#2}}}
\newcommand{\NormCzNd}[2]{\|#1\|_{\FsCz{#2}}}

\newcommand{\NormC}[3]{\BrN{#1}_{\FsC{#2}{#3}}}

\newcommand{\SNormC}[3]{\left|#1\right|_{\FsC{#2}{#3}}}

\newcommand{\indexSPH}{{\mbox{\rm\tiny SPH}}}
\newcommand{\ws}{w^\indexSPH}
\newcommand{\whs}{w_h^\indexSPH}

\newcommand{\IntpAppSPH}{\Pi_\h^{\indexSPH}}
\newcommand{\GradAppSPH}{\Grad_\h^{\indexSPH}}
\newcommand{\LapAppSPH}{\Delta_\h^{\indexSPH}}

\newcommand{\indexMPS}{{\mbox{\rm\tiny MPS}}}
\newcommand{\wm}{w^\indexMPS}
\newcommand{\whm}{w_h^\indexMPS}

\newcommand{\GradAppMPS}{\Grad_\h^{\indexMPS}}
\newcommand{\LapAppMPS}{\Delta_\h^{\indexMPS}}

\newcommand{\ParamMPSNZero}{\widehat{n}}
\newcommand{\ParamMPSLambda}{\widehat{\lambda}}

\newcommand{\mderivD}{{\rm D}}
\newcommand{\mderiv}[3]{
	\ifnum #1=1	\frac{\mderivD#3}{\mderivD{#2}}%
	\else{\frac{\mderivD^{#1}#3}{\mderivD^{#1}{#2}} %
	}\fi%
}

\newcommand{\StabFuncCher}{\Phi}
\newcommand{\StabFunc}[1]{\StabFuncCher}
\newcommand{\StabFuncMin}[1]{\StabFuncCher_{\min}}
\newcommand{\StabFuncMinApp}[1]{\StabFuncCher_{\rm app}}

\newcommand{\dN}{\mathbb{N}}

\newcommand{\dNz}{\mathbb{N}_0}

\newcommand{\dR}{\mathbb{R}}
\newcommand{\dRd}{\mathbb{R}^{\Dim}}

\newcommand{\dRp}{\mathbb{R}^{+}}
\newcommand{\dRpz}{\mathbb{R}^{+}_0}
\newcommand{\dZ}{\mathbb{Z}}

\newcommand{\xchar}{x}
\newcommand{\x}{\xchar}

\newcommand{\ychar}{y}
\newcommand{\y}{\ychar}

\newcommand{\zchar}{z}
\newcommand{\z}{\zchar}

\newcommand{\IntD}{\mathrm{d}}
\newcommand{\dx}{\IntD\x}
\newcommand{\dy}{\IntD\y}

\def\deq{\mathrel{\mathop:}=}%

\newcommand{\multiindexa}{\alpha}
\newcommand{\multiindexb}{\beta}

\newcommand{\NormDiscL}[3]{\Norm{#1}{\ell^{#2}(#3)}}

\newcommand{\DNOp}{\mathcal{D}}

\newcommand{\HDNMlOp}[2]{\DNOp}

\newcommand{\IndexSetXH}[2]{\Lambda_0(#1,#2)}
\newcommand{\IndexSetXHAst}[2]{\Lambda(#1,#2)}
\newtheorem{theorem}{Theorem}
\newtheorem{definition}[theorem]{Definition}
\newtheorem{lemma}[theorem]{Lemma}
\newtheorem{remark}[theorem]{Remark}
\newtheorem{corollary}[theorem]{Corollary}
\newtheorem{hypothesis}[theorem]{Hypothesis}
\newtheorem{proposition}[theorem]{Proposition}
\usepackage{mathtools}
\mathtoolsset{showonlyrefs=true}
\newcommand{\Space}{\quad}

\newcommand{\dMI}[1]{\mathbb{A}^{#1}}

\newcommand{\thmconst}{c}

\newcommand{\intA}{\ell}

\newcommand{\errorfuncA}[1]{I_{#1}}
\newcommand{\errorfuncB}[2]{I_{#1,#2}}
\newcommand{\errorfuncC}[1]{J_{#1}}
\newcommand{\prffuncA}[1]{E_{#1}}

\newcommand{\PvDmFunc}[2]{\xi^{\ast}_{#1}(#2)}
\newcommand{\NoiseChar}{\eta}
\newcommand{\Noise}[3]{\NoiseChar_{#1 #2}^{(#3)}}

\title[Truncation error estimates in a generalized particle method]{Truncation error estimates of approximate operators in a generalized particle method}

\author[Y. Imoto]{Yusuke Imoto${}^1$}
\address{${}^1$Institute for Advanced Study, Kyoto UniversityYoshida Konoe-cho, Sakyo-ku, Kyoto 6068501, Japan}
\email{imoto.yusuke.4e@kyoto-u.ac.jp}


\begin{document}
	
\begin{abstract}
	To facilitate the numerical analysis of particle methods, we derive truncation error estimates for the approximate operators in a generalized particle method. Here, a generalized particle method is defined as a meshfree numerical method that typically includes other conventional particle methods, such as smoothed particle hydrodynamics or moving particle semi-implicit methods. A new regularity of discrete parameters is proposed via two new indicators based on the Voronoi decomposition of the domain along with two hypotheses of reference weight functions. Then, truncation error estimates are derived for an interpolant, approximate gradient operator, and approximate Laplace operator in the generalized particle method. The convergence rates for these estimates are determined based on the frequency with which they appear in the regularity and hypotheses. Finally, the estimates are computed numerically and the results are shown to be in good agreement with the theoretical results. 
\end{abstract}
	
	\maketitle
	
\section{Introduction}
\label{sec:introduction}
Particle methods, such as the smoothed particle hydrodynamics (SPH) \cite{gingold1977smoothed,liu2010smoothed,lucy1977numerical} and moving particle semi-implicit (MPS) \cite{koshizuka1995particle,koshizuka1996moving,shakibaeinia2012mps} methods, are numerical methods for solving partial differential equations that are based on points called particles distributed in a domain. In such methods, an interpolant and several approximate differential operators are defined in terms of linear combinations of weighted interactions between neighboring particles. When such methods are applied to partial differential equations, the equations are effectively discretized in space. As the discretization procedure does not require mesh generation in the domain, particle methods can be applied to moving boundary problems, such as the deformation and destruction of structures \cite{benz1995simulations,monaghan1991simulation} and flow problems associated with free surfaces \cite{monaghan1994simulating,murotani2014development}. 

The accuracy of particle methods has been widely researched. 
From an engineering perspective, many studies have been conducted into the convergence of such methods in practical applications, such as Amicarelli \cite{amicarelli2011sphintp,amicarelli2011sphgrad}, Fulk \cite{fulk1996analysis}, and Quinlan et al. \cite{quinlan2006truncation}. 
On the other hand, few studies in the literature have presented numerical analyses of these methods from a mathematical perspective. 
In the 1980s, Mas-Gallic and Ravirt \cite{mas1987particle} and Ravirt \cite{raviart1985analysis} provided error estimates for particle methods when applied to parabolic and hyperbolic partial differential equations on unbounded domains. In the 2000s, Moussa and Via \cite{moussa2000convergence} and Moussa \cite{moussa2006convergence} provided error estimates of nonlinear conservation lows on bounded domains. In their work, the time integrations of the particle positions and volumes were obtained by solving the differential equations with respect to advection fields. However, as their method is only applicable to problems described by solvable differential equations, it cannot be used with other problems, such as those involving the Navier--Stokes equations. 

Sometime later, Ishijima and Kimura \cite{ishijima2010mashfreeE} developed a truncation error estimate for an approximate gradient operator in MPS. By introducing a regularity for particle distributions based on an indicator called the equivolume partition radius, they determined the conditions that depend solely on the space distributions of the particles. However, a practical limitation is that the indicator cannot be computed. 

In previous works, we established truncation error estimates for an interpolant, approximate gradient operator, and approximate Laplace operator of a generalized particle method in which the particle volumes were given as Voronoi volumes \cite{imoto2016teintpV,imoto2017tedifV}. A generalized particle method is a numerical method that typically includes conventional particle methods, such as SPH and MPS. In previous studies, we derived truncation error estimates by introducing a regularity using an indicator known as the covering radius, which is used in the numerical analysis of meshfree methods based on moving least-square methods and radial basis functions \cite{levin1998errorMLS,schaback1995error,wendland2005scattered}. 
Although the formulations and conditions in those works are computable, they are difficult to deploy in practical computations as the computational costs associated with particle volumes based on Voronoi decomposition are high. 

The focus of the current work was to analyze particle methods under more practical conditions by extending our results to cases with commonly used particle volumes. 
We also introduce another indicator of particle volumes, which we refer to as a Voronoi deviation, that represents the deviation between particle volumes and Voronoi volumes. 
Then, utilizing the Voronoi deviation, we extend the regularity and introduce two hypotheses of reference weight functions. 
Using the regularity and hypotheses, we derive truncation error estimates of the interpolant, approximate gradient operator, and approximate gradient operator of the generalized particle method. 
Finally, we numerically analyze our estimates and compare the results to those from the theory. 

The remainder of this paper is organized as follows.
The interpolant and approximate operators of the generalized particle method are introduced in Section \ref{sec:formulation}. A regularity describing the family of discrete parameters is discussed in Section \ref{sec:truncation_error_estimate}, after which we propose our primary theorem with respect to the truncation error estimates and provide some corollaries. Then, the primary theorem is proven in Section \ref{sec:Proof}, numerical results are detailed in Section \ref{sec:Numerical_results}, and some concluding remarks are outlined in Section \ref{sec:conclusions}. 

In the remainder of this section, we describe some notation and define some relevant function spaces. 
Let $\dRp$, $\dRpz$, and $\dNz$ be the set of positive real numbers, the set of nonnegative real numbers, and the set of nonnegative integers, respectively. 
Let $\Dim$ be the dimension of a space. 
Let $\dMI{\Dim}$ be the set of all $\Dim$-dimensional multi-indices. 
If there is no ambiguity, the symbol $|\cdot|$ is used to denote the following: $|x|$ denotes the Euclidean norm for a vector $x\in\dRd$; $|S|$ denotes the volume of $S$ for open set $S\subset\dRd$; $|\alpha|$ denotes $|\alpha|\deq\alpha_1+\alpha_2+\dots+\alpha_\Dim$ for $\alpha\in\dMI{\Dim}$. 
For $S\subset\dRd$, let $\diam(S)$ be $\diam(S)\deq \sup\Set{|\x-\y|}{\x,\y\in S}$. 
For $S\subset\dRd$, let $\FsCz{\Ol{S}}$ be the space of real continuous functions defined in $\Ol{S}$ with the norm $\NormCz{\,\cdot\,}{\Ol{S}}$ defined as
\begin{align}
	\NormCz{\f}{\Ol{S}} &\deq \max_{\x\in\Ol{S}}\BrA{\f(\x)}.
\end{align}
For $S\subset\dRd$ and $\ell\in\dN$, let $\FsC{\ell}{\Ol{S}}$ be the space of functions in $\FsCz{\Ol{S}}$ with derivatives up to the $\ell$th order with its seminorm $\SNormC{\,\cdot\,}{\ell}{\Ol{S}}$ and norm $\NormC{\,\cdot\,}{\ell}{\Ol{S}}$ defined as
\begin{align}
	\BrA{\f}_{\FsC{\ell}{\Ol{S}}} &\deq \max_{\alpha\in\dMI{\Dim}, |\alpha|=\ell}\BrN{D^\alpha \f}_{\FsCz{\Ol{S}}},
	\\
	\BrN{\f}_{\FsC{\ell}{\Ol{S}}} &\deq \max_{j=0,1,\dots,\ell} \SNormC{\f}{j}{\Ol{S}},
\end{align}
respectively.
Here $D^\alpha\f\deq\partial_1^{\alpha_1} \partial_2^{\alpha_2} \dots \partial_\Dim^{\alpha_\Dim}\f$ with multi-index $\alpha=(\alpha_1,\alpha_2,\dots,\alpha_\Dim)$. 

\section{Approximate operators in a generalized particle method}
\label{sec:formulation}
\begin{figure}[t]
	\includegraphics[width=55mm,bb=13 17 1837 1185]{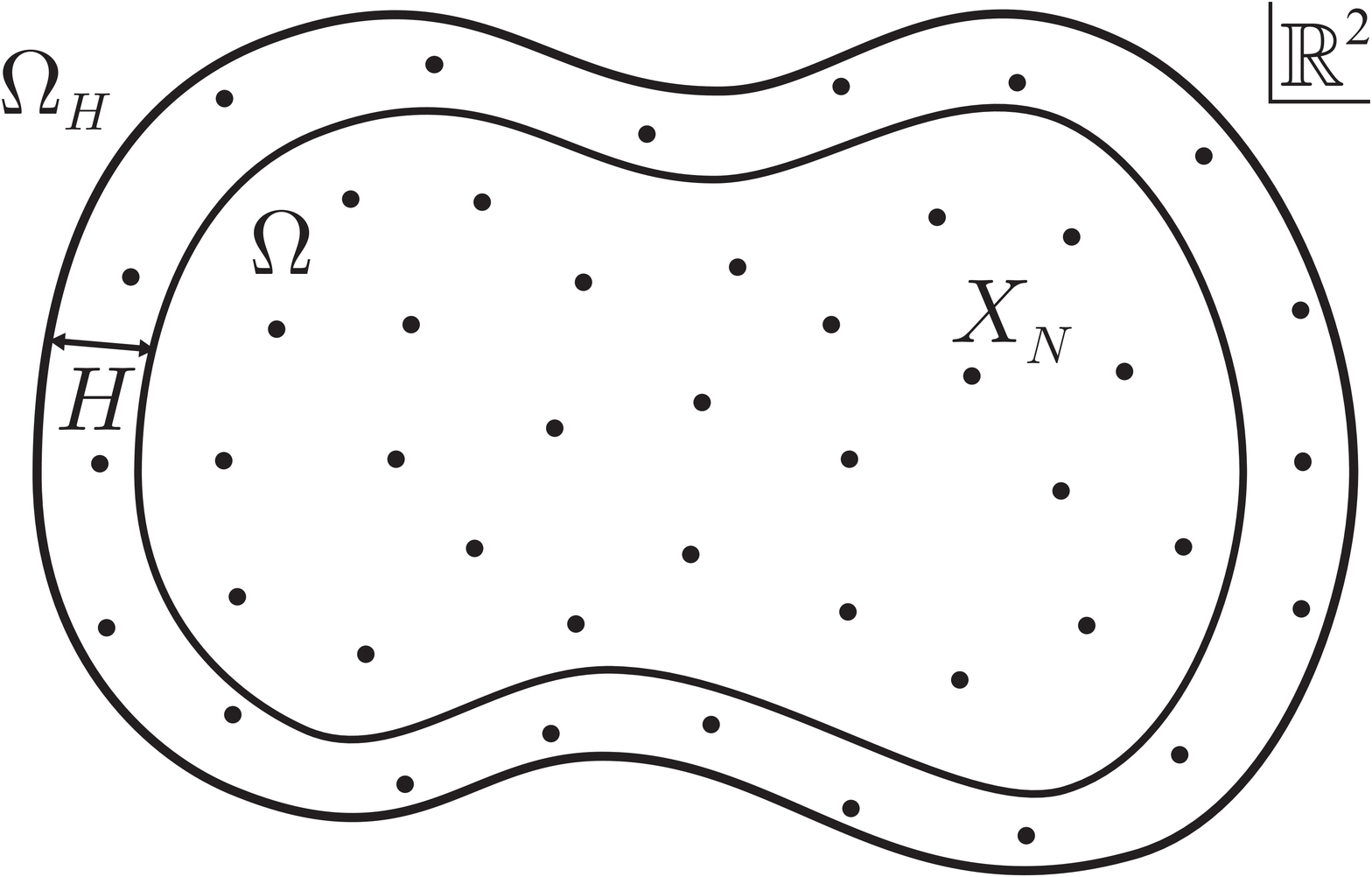}
	\caption{Particle distribution $\PtSet$ in $\DmH\,(\subset\dR^2)$.}
	\label{fig_particle_distribution}
\end{figure}

Let $\Dm$ be a bounded domain in $\dRd$. 
Let $\H$ be a fixed positive number. 
For $\Dm$ and $\H$, we define extended domain $\DmH$ as 
\begin{equation}
	\DmH \deq \left\{\x\in\dRd \midd \ex \y\in\Dm \st |\x-\y|<\H\right\}. 
\end{equation}
For $\N\in\dN$, we define a particle distribution $\PtSet$ and particle volume set $\PvSet$ as
\begin{equation}
	\PtSet \deq \Set{\Pt{\i}\in\DmH}{\i=1,2,\dots,\N,\quad\Pt{\i}\neq\Pt{\j}\,(\i\neq\j)}, 
	\label{def:PtSet}
\end{equation}
\begin{equation}
	\PvSet \deq \Set{\Pv{\i}\in\dRp}{\i=1,2,\dots,\N,\quad\sum_{\i=1}^\N\Pv{\i}=\Meas{\DmH}}, 
	\label{def:PvSet}
\end{equation}
respectively. 
We refer to $\Pt{\i}\in\PtSet$ and $\Pv{\i}\in\PvSet$ as a particle and particle volume, respectively.   
An example of the particle distribution $\PtSet$ in $\DmH\,(\subset\dR^2)$ is shown in Figure \ref{fig_particle_distribution}. 

We define an admissible reference weight function set $\FsWeightFunc$ as
\begin{equation}
	\FsWeightFunc \deq \Set{\w\in\FsCz{\dRpz}}{\supp(\w)=[0,1],\,\int_{\dRd}\w(|\x|)\dx=1, \mbox{absolutely continuous}}, 
\end{equation}
we refer to $\w\in\FsWeightFunc$ as a reference weight function, and 
we define the influence radius $\hN\in\dR$ as satisfying $0<\hN<\H$ and $\hN\ra0\,(\N\ra\infty)$. 
If there is no ambiguity, we denote $\hN$ as $\h$. 
For reference weight function $\w$ and influence radius $\h$, we define the weight function $\wh\in\FsCz{\dRpz}$ as
\begin{equation}
	\wh(r) \deq \frac{1}{\h^{\Dim}}\w\left(\frac{r}{\h}\right). 
	\label{def:wh}
\end{equation}
Note that the weight function $\wh$ satisfies
\begin{align}
	\supp(\wh) = [0,\h], \qquad \int_{\dRd}\wh(|\x|)\dx=1, 
\end{align}
and is absolutely continuous. 

For $\f\in\FsCz{\DmHOl}$, we define interpolant $\IntpApp{}$, approximate gradient operator $\GradApp{}$, and approximate Laplace operator $\LapApp{}$ as
\begin{align}
	\IntpApp{} \f (\x) &\deq \sum_{\i\in\IndexSetXH{\x}{\h}}\Pv{\i} \f(\Pt{\i}) \wh(|\Pt{\i}-\x|),
	\label{ih_def}\\
	\GradApp{} \f (\x) &\deq \Dim  \sum_{\i\in\IndexSetXHAst{\x}{\h}}\Pv{\i} \frac{\f(\Pt{\i})-\f(\x)}{|\Pt{\i}-\x|}\frac{\Pt{\i}-\x}{|\Pt{\i}-\x|} \wh(|\Pt{\i}-\x|),
	\label{gh_def}\\
	\LapApp{} \f (\x) &\deq 2\Dim\sum_{\i\in\IndexSetXHAst{\x}{\h}}\Pv{\i} \frac{\f(\Pt{\i})-\f(\x)}{|\Pt{\i}-\x|^2} \wh(|\Pt{\i}-\x|), 
	\label{lh_def}
\end{align}
respectively. 
Here, for $\x\subset\dRd$ and $r\in\dRp$, $\IndexSetXH{\x}{r}$ and $\IndexSetXHAst{\x}{r}$ are index sets of particles defined as
\begin{align}
	\IndexSetXH{\x}{r} &\deq \Set{\i=1,2,\dots,\N}{0\leq|\x-\Pt{\i}|<r}, \\
	\IndexSetXHAst{\x}{r} &\deq \Set{\i=1,2,\dots,\N}{0<|\x-\Pt{\i}|<r}, 
\end{align}
respectively. 

As discussed later in Appendix \ref{sec:appendix_conventional_particle_methods}, the approximate operators \eqref{ih_def}, \eqref{gh_def}, and \eqref{lh_def} indicate a wider class of approximate operators of particle methods that those in the SPH and MPS methods. 
Therefore, we refer to the approximate operators \eqref{ih_def}, \eqref{gh_def}, and \eqref{lh_def} as generalized approximate operators and to a particle method that uses them as a generalized particle method. 

\section{Truncation error estimates of approximate operators}
\label{sec:truncation_error_estimate}
\begin{figure}[t]
	\includegraphics[width=55mm,bb=0 0 653.7mm 420.9mm]{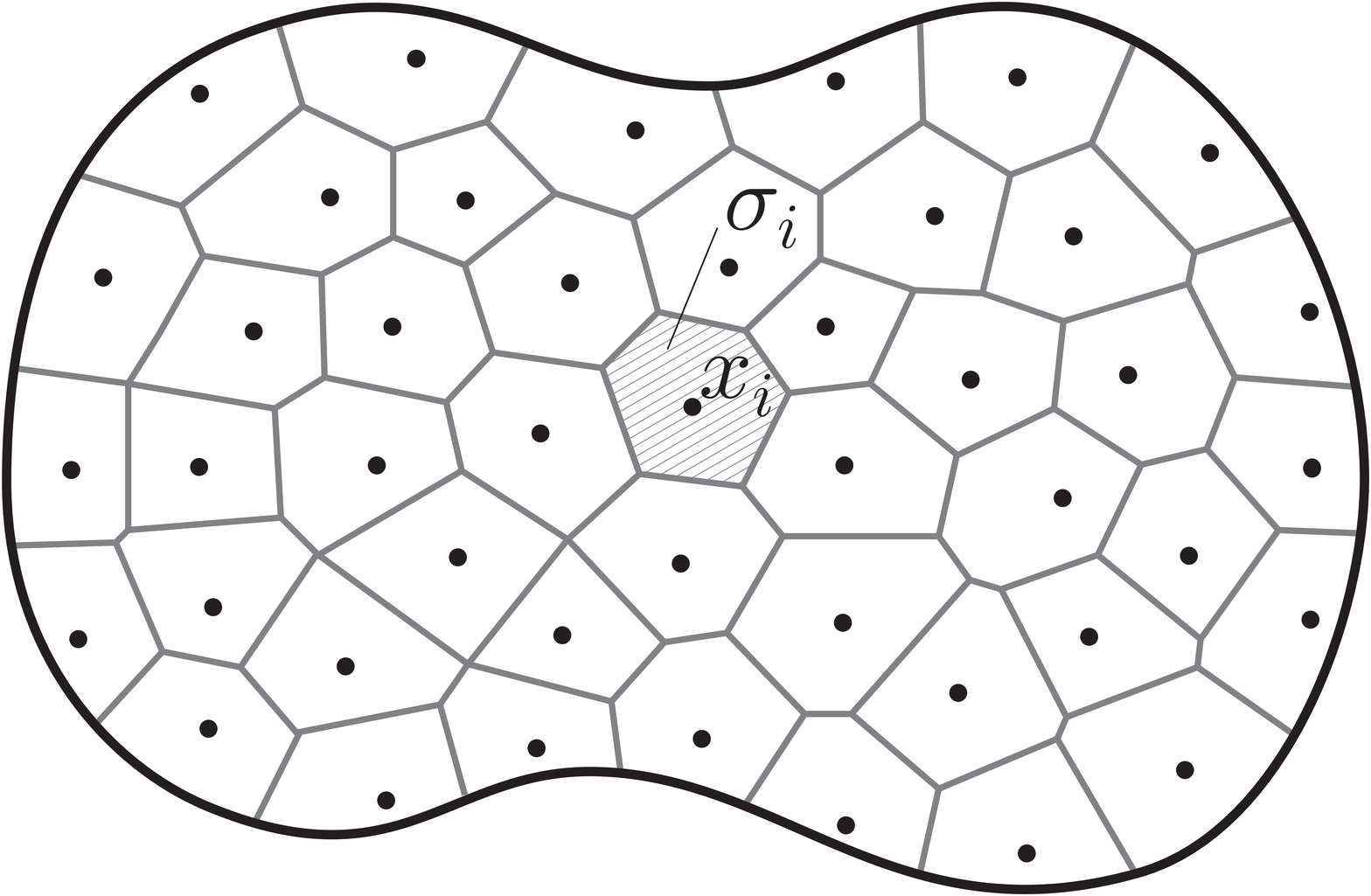}
	\caption{Example of the Voronoi decomposition of $\DmH$ associated with the particle distribution $\PtSet$. }
	\label{fig:Voronoi_decomposition}
\end{figure}

We first introduce a regularity of discrete parameters. 
Let $\{\VoroCell{\i}\}$ be the Voronoi decomposition of $\DmH$ associated with the particle distribution $\PtSet$, where $\VoroCell{\i}$ is the Voronoi region defined as 
\begin{align}
	\VoroCell{\i} \deq \Set{\x\in\DmH}{|\Pt{\i}-\x|<|\Pt{\j}-\x|,~\fa \Pt{\j}\in\PtSet\,(\j\neq \i) },\quad \i=1,2,\dots,\N. 
	\label{def_Voronoi_cell}
\end{align}
We define a particle volume decomposition $\PvDmSet$ as a decomposition of $\DmH$ associated with particle volume set $\PvSet$ defined as
\begin{equation}
	\PvDmSet\deq\Set{\PvDm{\i}\subset\DmH}{\Meas{\PvDm{\i}}= \Pv{\i},\quad\bigcup_{\i=1}^\N\PvDmOl{\i}=\DmHOl,\quad\PvDm{\i}\cap\PvDm{\j}=\emptyset\,(\i \neq \j)}. 
\end{equation}
An example of the Voronoi decomposition of $\DmH$ associated with the particle distribution $\PtSet$ is shown in Figure \ref{fig:Voronoi_decomposition}.  
We define a covering radius $\IndPtSet$ for particle distribution $\PtSet$ as
\begin{align}
	\IndPtSet \deq \max_{\i=1,2,\dots,\N} \sup_{\x\in\VoroCell{\i}}|\Pt{\i}-\x|. 
	\label{def:covering_radius}
\end{align}
Moreover, we define a Voronoi deviation $\IndPvSet$ for the particle distribution $\PtSet$ and the particle volume set $\PvSet$ as
\begin{align}
	\IndPvSet\deq \inf_{\PvDmSet} \IndPvSetLocal{\PvDmSet}
	\label{def:Voronoi_deviation}
\end{align}
with
\begin{equation}
	\IndPvSetLocal{\PvDmSet} \deq \maxi{\i}\BrM{\sum_{\j=1}^\N\dfrac{\Meas{\VoroCell{\i}\cap\PvDm{\j}}+\Meas{\PvDm{\i}\cap\VoroCell{\j}}}{\Meas{\VoroCell{\i}}}|\Pt{\i}-\Pt{\j}|}. 
	\label{def:Voronoi_deviation_local}
\end{equation}
Then, we define a regularity for a family consisting of a particle distribution $\PtSet$, particle volume set $\PvSet$, and influence radius $\h$ as follows: 
\begin{definition}
	\label{def:regular}
	Family $\{(\PtSet, \PvSet, \hN)\}_{\N\ra\infty}$ is said to be regular with order $\RegOrder\,(\RegOrder\geq 1)$ if there exists a positive constant $\regularconst$ such that
	\begin{align}
		\hN^{\RegOrder} \geq \regularconst (\IndPtSet+\IndPvSet),\qquad \fa\N\in\dN. 
		\label{def:regular:cond}
	\end{align}
\end{definition}

\begin{figure}[t]
	\includegraphics[width=120mm,bb=0 0 1551mm 559mm]{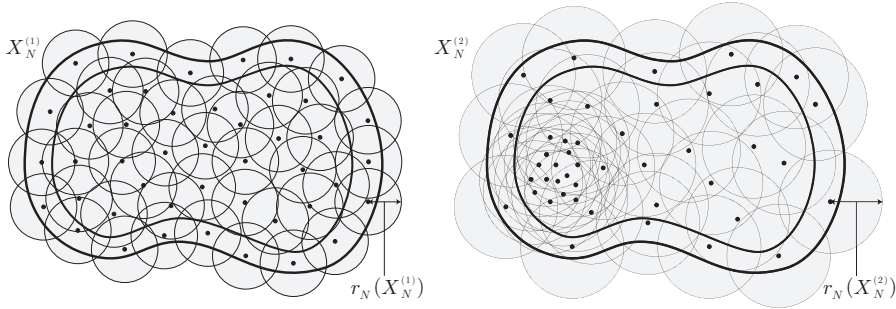}
	\caption{Two examples of covering radii $\IndPtSet$ for particle distributions with same number of particles. 
		The covering radius $\IndPtSet$ for the uniform particle distribution \rmn{(}left\rmn{)} is smaller than that for the uneven particle distribution \rmn{(}right\rmn{)}. }
	\label{fig:covering_radius}
\end{figure}

\begin{remark}
	As shown in Figure \ref{fig:covering_radius}, the covering radius $\IndPtSet$ becomes large in the case of a particle distribution with both dense and sparse regions. 
	Therefore, the covering radius $\IndPtSet$ can be considered as an indicator representing the uniformness of particle distribution $\PtSet$. 
\end{remark}

\begin{remark}
	A Voronoi deviation $\IndPvSet$ equals zero if and only if the particle volumes are given as the Voronoi volume ($\Pv{\i}=\Meas{\VoroCell{\i}}$). 
	Moreover, the Voronoi deviation $\IndPvSet$ becomes large if the particle volumes are given as values far from the Voronoi volumes. 
	Therefore, the Voronoi deviation $\IndPvSet$ can be regarded as an indicator of the deviation between the particle volume set and the Voronoi volume set. 
\end{remark}

\begin{remark}
	For a given family $\{(\PtSet, \PvSet, \hN)\}_{\N\ra\infty}$ and given constant $\RegOrder$ $(\RegOrder\geq1)$, it is possible to determine whether or not the family is regular with order $\RegOrder$ as the covering radius $\IndPtSet$ and Voronoi deviation $\IndPvSet$ are absolutely computable, as shown in Appendix \rmn{\ref{sec:appendix_compt_procedure_indicator}}. 
\end{remark}

Next, we introduce two hypotheses of reference weight function $\w$:  
\begin{hypothesis}
	\label{hypo:ref_weight:order}
	For $\WeightOrder\in\dN$, the reference weight function $\w$ satisfies for all $\multiindexa\in\dMI{\Dim}$ with $1\leq|\multiindexa|\leq \WeightOrder$, 
	\begin{align}
		\int_{\dRd} \x^\multiindexa \w(|\x|)\dx = 0. 
		\label{w_hypo_order_eq_01}
	\end{align}
\end{hypothesis}
\begin{hypothesis}
	\label{hypo:ref_weight:smooth}
	For $\WeightOrderSmooth\in\dNz$, the reference weight function $\w$ satisfies
	\begin{align}
		\max\BrM{\sup_{r\in(0,1)}\BrANd{\wn{\WeightOrderSmooth+1}(r)}, \sup_{r\in(0,1)}\BrA{\dfrac{d}{dr}\wn{\WeightOrderSmooth}(r)}} < \infty,
		\label{w_hypo_smooth_01_eq}
	\end{align}
	where for $j\in\dNz$, $\wn{j}(r): (0,\infty)\ra\dR$ is defined as
	\begin{align}
		\wn{j}(r) \deq
		\begin{cases}
			\ds \lim_{s\downarrow 0}\dfrac{\w(s)}{s^{j}},\qquad & r=0, 
			\\
			\dfrac{\w(r)}{r^{j}},\qquad & r>0. 
		\end{cases}
		\label{def:wn}
	\end{align}
\end{hypothesis}

\begin{remark}
	All reference functions $\w\in\FsWeightFunc$ satisfy Hypothesis \ref{hypo:ref_weight:order} with $\WeightOrder=1$. 
	Moreover, for all $\WeightOrder\in\dN$ and $\WeightOrderSmooth\in\dN$,  reference weight functions satisfying Hypothesis \ref{hypo:ref_weight:order} with $\WeightOrder$ and Hypothesis \ref{hypo:ref_weight:smooth} with $\WeightOrderSmooth$ can be constructed as shown in Appendix \ref{sec:appendix_const_ref_weight_func}. 
\end{remark}

We now state a theorem that defines truncation error estimates of approximate operators in the generalized particle method with a continuous norm: 
\begin{theorem}
	\label{thm:terror_cont_norm}
	Suppose that family $\{(\PtSet,\PvSet,\hN)\}_{\N\ra\infty}$ is regular with order $\RegOrder\,(\RegOrder\geq1)$ and that reference weight function $\w$ satisfies Hypothesis \ref{hypo:ref_weight:order} with $\WeightOrder$. 
	Then, there exists a positive constant $\thmconst$ independent of $\N$ such that
	\begin{align}
		\NormCz{\f-\IntpApp{} \f}{\DmOl} & \leq \thmconst\, \h^{\min\{\RegOrder-1,\WeightOrder+1\}} \NormC{\f}{\WeightOrder+1}{\DmHOl},\quad \f\in \FsC{\WeightOrder+1}{\DmHOl}. 
		\label{ineq:terror_cont_norm_intp}
	\end{align}
	In addition, if $\w\in\FsWeightFunc$ satisfies Hypothesis \ref{hypo:ref_weight:smooth} with $\WeightOrderSmooth=0$, then we have 
	\begin{align}
		\NormCz{\Grad\f-\GradApp{}\f}{\DmOl} \leq \thmconst\, \h^{\min\{\RegOrder-1,\WeightOrder+1\}} \NormC{\f}{\WeightOrder+2}{\DmHOl},\quad \f\in \FsC{\WeightOrder+2}{\DmHOl}, 
		\label{ineq:terror_cont_norm_grad}
	\end{align}
	and if $\w\in\FsWeightFunc$ satisfies Hypothesis \ref{hypo:ref_weight:smooth} with $\WeightOrderSmooth=1$, then we have 
	\begin{align}
		\NormCz{\Lap\f-\LapApp{}\f}{\DmOl} \leq \thmconst\, \h^{\min\{\RegOrder-2,\WeightOrder+1\}} \BrN{\f}_{\FsC{\WeightOrder+3}{\DmHOl}},\quad \f\in \FsC{\WeightOrder+3}{\DmHOl}. 
		\label{ineq:terror_cont_norm_lap}
	\end{align}
\end{theorem}

The proof of Theorem \ref{thm:terror_cont_norm} is presented in the next section. 
As shown in the corollaries in Appendix \ref{sec:appendix_conventional_particle_methods}, the approximate operators commonly used in the SPH and MPS methods are valid for Theorem \ref{thm:terror_cont_norm} under appropriate settings. 


\section{Proof of truncation error estimates}
\label{sec:Proof}
The following notation will be used in the subsequent proof of Theorem \ref{thm:terror_cont_norm}. 
Hereafter, let $\gconst$ be a generic positive constant independent of $\N$. 
For $\multiindexa\in\dMI{\Dim}$, set $\errorfuncA{\multiindexa}$ as 
\begin{equation}
	\errorfuncA{\multiindexa}(\x) \deq \sum_{\i\in\IndexSetXH{\x}{\h}}\Pv{\i} (\Pt{\i}-\x)^{\multiindexa} \wh(|\Pt{\i}-\x|) - \int_{\dRd}\y^\multiindexa \wh(|\y|)\dy,\qquad \x\in\DmOl. 
\end{equation}
For $\multiindexa\in\dMI{\Dim}$ and $\intA\in\dN$, set $\errorfuncB{\multiindexa}{\intA}$ as
\begin{equation}
	\errorfuncB{\multiindexa}{\intA}(\x) \deq \sum_{\i\in\IndexSetXHAst{\x}{\h}} \Pv{\i} \frac{(\Pt{\i}-\x)^{\multiindexa}}{|\Pt{\i}-\x|^{\intA}} \wh(|\Pt{\i}-\x|) - \int_{\dRd} \frac{\y^\multiindexa}{|\y|^\intA} \wh(|\y|)\dy,\quad  \x\in\DmOl. 
\end{equation}
For $\intA\in\dN$, set $\errorfuncC{\intA}$ as
\begin{equation}
	\errorfuncC{\intA}(\x) \deq \sum_{\i\in\IndexSetXH{\x}{\h}} \Pv{\i} |\Pt{\i}-\x|^{\intA} |\wh(|\Pt{\i}-\x|)|,\qquad\x\in\DmOl. 
\end{equation}
We now present the following lemma. 

\begin{lemma}
	\label{lem:terror_cont_norm}
	Suppose that $\w\in\FsWeightFunc$ satisfies Hypothesis \ref{hypo:ref_weight:order} with $\WeightOrder$. 
	Then, there exists a positive constant $\thmconst$ independent of $\N$ such that 
	\begin{equation}
		\NormCz{\f-\IntpApp{} \f}{\DmOl}
		\leq \gconst\BrS{\sum_{1\leq|\multiindexa|\leq\WeightOrder}\NormCz{\errorfuncA{\multiindexa}}{\DmOl}+ \NormCz{\errorfuncC{\WeightOrder+1}}{\DmOl}}\NormC{\f}{\WeightOrder+1}{\DmHOl},\quad\f\in\FsC{\WeightOrder+1}{\DmHOl},
		\label{lem:terror_cont_norm:intp}
	\end{equation}
	\begin{equation}
		\NormCz{\Grad\f-\GradApp{}\f}{\DmOl}
		\leq\gconst\BrS{\sum_{2\leq |\multiindexa|\leq\WeightOrder+2} \NormCz{\errorfuncB{\multiindexa}{2}}{\DmOl}+ \NormCz{\errorfuncC{\WeightOrder+1}}{\DmOl}}\NormC{\f}{\WeightOrder+2}{\DmHOl},\quad\f\in\FsC{\WeightOrder+2}{\DmHOl},
		\label{lem:terror_cont_norm:grad}
	\end{equation}
	\begin{equation}
		\NormCz{\Lap\f-\LapApp{}\f}{\DmOl}
		\leq \gconst\BrS{\sum_{1\leq |\multiindexa|\leq\WeightOrder+3} \NormCz{\errorfuncB{\multiindexa}{2}}{\DmOl}+ \NormCz{\errorfuncC{\WeightOrder+1}}{\DmOl}}\NormC{\f}{\WeightOrder+3}{\DmHOl},\quad\f\in\FsC{\WeightOrder+3}{\DmHOl}.
		\label{lem:terror_cont_norm:lap}
	\end{equation}
\end{lemma}

\begin{proof}
	First, we prove \eqref{lem:terror_cont_norm:intp}. 
	We fix $\x\in\DmOl$. 
	Then, let $\OBall{\x}{r}$ be the open ball in $\dRd$ with center $\x$ and radius $r$, i.e.,
	\begin{equation}
		\OBall{\x}{r}\deq \Set{\y\in\dRd}{|\x-\y|<r}. 
	\end{equation}
	From $\h<\H$, we have $\OBall{\x}{\h}\subset\DmH$.
	Then, for all $\f\in\FsC{\ell+1}{\DmHOl}$ $(\ell\in\dN)$ and $\y\in\OBall{\x}{\h}$, we obtain the Taylor expansion of $\f$ as
	\begin{align}
		& \f(\Pt{\i}) = \sum_{0\leq |\multiindexa|\leq\ell} \frac{D^{\multiindexa}\f(\x)}{\multiindexa!}(\x-\Pt{\i})^{\multiindexa} + \sum_{|\multiindexa|=\WeightOrder+1}(\x-\Pt{\i})^{\multiindexa} R_{\multiindexa}(\Pt{\i},\x), 
		\label{taylor01}\\
		& R_{\multiindexa}(\Pt{\i},\x) \deq \frac{|\multiindexa|}{\multiindexa !} \int_0^1 (1-t)^{|\multiindexa|-1}D^{\multiindexa}\f(t\Pt{\i}+(1-t)\x)dt. 
		\label{regidual01}
	\end{align}
	From \eqref{ih_def} and \eqref{taylor01} with $\ell=\WeightOrder+1$, we have 
	\begin{align}
		\IntpApp{} \f(\x) 
		& = \sum_{0\leq |\multiindexa|\leq\WeightOrder} \frac{D^{\multiindexa}\f(\x)}{\multiindexa!} \sum_{\i\in\IndexSetXH{\x}{\h}} \Pv{\i} (\x-\Pt{\i})^{\multiindexa} \wh(|\x-\Pt{\i}|)
		\nn\\
		&\quad+\sum_{|\multiindexa|=\WeightOrder+1}\sum_{\i\in\IndexSetXH{\x}{\h}}R_{\multiindexa}(\Pt{\i},\x)\Pv{\i} (\x-\Pt{\i})^{\multiindexa} \wh(|\x-\Pt{\i}|). 
	\end{align}
	Moreover, by Hypothesis \ref{hypo:ref_weight:order}, we have 
	\begin{align}
		\IntpApp{} \f(\x) - \f(\x)
		& = \sum_{0\leq |\multiindexa|\leq\WeightOrder} \frac{D^{\multiindexa}\f(\x)}{\multiindexa!} \errorfuncA{\multiindexa}(\x)
		\nn\\
		&\quad+\sum_{|\multiindexa|=\WeightOrder+1}\sum_{\i\in\IndexSetXH{\x}{\h}}R_{\multiindexa}(\Pt{\i},\x)\Pv{\i} (\x-\Pt{\i})^{\multiindexa} \wh(|\x-\Pt{\i}|). 
		\label{prf:lem:terror_cont_norm:intp01}
	\end{align}
	Because 
	\begin{align}
		|R_{\multiindexa}(\y,\z)| \leq \frac{1}{\multiindexa!}\SNormC{\f}{|\multiindexa|}{\DmHOl}, \qquad  \y\in\DmOl,~\z\in\OBall{\y}{\h},~\multiindexa\in\dMI{\Dim},
		\label{ineq:R_a[f](y_x)}
	\end{align}
	we have
	\begin{equation}
		\BrA{\sum_{|\multiindexa|=\WeightOrder+1}\sum_{\i\in\IndexSetXH{\x}{\h}}R_{\multiindexa}(\Pt{\i},\x)\Pv{\i} (\x-\Pt{\i})^{\multiindexa} \wh(|\x-\Pt{\i}|)}
		\leq\gconst|\errorfuncC{\WeightOrder+1}(\x)|\SNormC{\f}{\WeightOrder+1}{\DmHOl}.
		\label{eq:int_y^a_wh(|y|)dy}
	\end{equation}
	Moreover, we have
	\begin{align}
		\BrA{\sum_{0\leq |\multiindexa|\leq\WeightOrder} \frac{D^{\multiindexa}\f(\x)}{\multiindexa!} \errorfuncA{\multiindexa}(\x)}
		&\leq\gconst\NormC{\f}{\WeightOrder}{\DmOl}\sum_{1\leq|\multiindexa|\leq\WeightOrder}\BrA{\errorfuncA{\multiindexa}(\x)}. 
		\label{ineq:Ja_intp}
	\end{align}
	Therefore, from \eqref{prf:lem:terror_cont_norm:intp01}, \eqref{eq:int_y^a_wh(|y|)dy}, and \eqref{ineq:Ja_intp}, we obtain \eqref{lem:terror_cont_norm:intp}. 
	
	Next, we prove \eqref{lem:terror_cont_norm:grad}.	
	From \eqref{gh_def} and \eqref{taylor01} with $\ell=\WeightOrder+2$, we have 
	\begin{align}
		\GradApp{} \f(\x) 
		&=\Dim\sum_{1\leq|\multiindexa|\leq\WeightOrder+1}\frac{D^{\multiindexa}\f(\x)}{\multiindexa!}\sum_{\i\in\IndexSetXHAst{\x}{\h}}\Pv{\i}\frac{(\x-\Pt{\i})(\x-\Pt{\i})^{\multiindexa}}{|\x-\Pt{\i}|^2} \wh(|\x-\Pt{\i}|)
		\nn\\
		&\quad+\Dim\sum_{|\multiindexa|=\WeightOrder+2}\sum_{\i\in\IndexSetXHAst{\x}{\h}}R_{\multiindexa}(\Pt{\i},\x)\Pv{\i}\frac{(\x-\Pt{\i})(\x-\Pt{\i})^{\multiindexa}}{|\x-\Pt{\i}|^2}\wh(|\x-\Pt{\i}|). 
	\end{align}
	Because for $\multiindexb\in\dMI{\Dim}$ with $|\multiindexb|=2$, 
	\begin{align}
		\Dim \int_{\dRd} \frac{\y^{\multiindexb}}{|\y|^2}\wh(|\y|)\dy =
		\begin{cases}
			1,\quad & \mbox{all elements of $\multiindexb$ are even},
			\\
			0,\quad & \mbox{otherwise}, 
		\end{cases}
		\label{eq:d_int_y^b_y^-2dy}
	\end{align}
	we have
	\begin{align}
		\Dim\sum_{|\multiindexa|=1}\frac{D^{\multiindexa}\f(\x)}{\multiindexa!}\int_{\dRd} \frac{y\y^\multiindexa}{|\y|^2} \wh(|\y|)\dy = \Grad\f(\x). 
		\label{eq:grad_rhs1}
	\end{align}
	Hypothesis \ref{hypo:ref_weight:order} with $\WeightOrder$ yields
	\begin{align}
		\int_{\dRd} \frac{\y \y^\multiindexa}{|\y|^2} \wh(|\y|)\dy = 0\qquad\mbox{$\multiindexa\in\dMI{\Dim}$ with $2\leq|\multiindexa|\leq\WeightOrder+1$}. 
		\label{eq:int_y^a_y^-2_wh_dy=0}
	\end{align}
	From \eqref{eq:grad_rhs1} and \eqref{eq:int_y^a_y^-2_wh_dy=0}, we have 
	\begin{align}
		\GradApp{}\f(\x)-\Grad\f(\x)&=-\Dim\sum_{1\leq|\multiindexa|\leq\WeightOrder+1}\frac{D^{\multiindexa}\f(\x)}{\multiindexa!}\int_{\dRd} \frac{\y \y^\multiindexa}{|\y|^2} \wh(|\y|)\dy
		\\
		&\quad+\Dim\sum_{1\leq|\multiindexa|\leq\WeightOrder+1}\frac{D^{\multiindexa}\f(\x)}{\multiindexa!}\sum_{\i\in\IndexSetXHAst{\x}{\h}}\Pv{\i}\frac{(\x-\Pt{\i})(\x-\Pt{\i})^{\multiindexa}}{|\x-\Pt{\i}|^2} \wh(|\x-\Pt{\i}|)
		\\
		&\quad+\Dim\sum_{|\multiindexa|=\WeightOrder+2}\sum_{\i\in\IndexSetXHAst{\x}{\h}}R_{\multiindexa}(\Pt{\i},\x) \Pv{\i}\frac{(\x-\Pt{\i})(\x-\Pt{\i})^{\multiindexa}}{|\x-\Pt{\i}|^2}\wh(|\x-\Pt{\i}|). 
		\label{prf:lem:terror_cont_norm:grad01}
	\end{align}
	From \eqref{ineq:R_a[f](y_x)}, we have
	\begin{equation}
		\BrA{\sum_{|\multiindexa|=\WeightOrder+2}\sum_{\i\in\IndexSetXHAst{\x}{\h}}R_{\multiindexa}(\Pt{\i},\x)\Pv{\i}\frac{(\x-\Pt{\i})(\x-\Pt{\i})^{\multiindexa}}{|\x-\Pt{\i}|^2} \wh(|\x-\Pt{\i}|)}
		\leq\gconst|\errorfuncC{\WeightOrder+1}(\x)|\SNormC{\f}{\WeightOrder+2}{\DmHOl}. 
		\label{ineq:R_grad}
	\end{equation}
	Moreover, we have
	\begin{equation}
		\sum_{1\leq|\multiindexa|\leq\WeightOrder+1}\BrA{\sum_{\i\in\IndexSetXHAst{\x}{\h}}\Pv{\i}\frac{(\x-\Pt{\i})(\x-\Pt{\i})^{\multiindexa}}{|\x-\Pt{\i}|^2}\wh(|\x-\Pt{\i}|)-\int_{\dRd}\frac{\y\y^\multiindexa}{|\y|^2}\wh(|\y|)\dy}\leq\gconst\sum_{2\leq |\multiindexa|\leq\WeightOrder+2} \BrANd{\errorfuncB{\multiindexa}{2}(\x)}. 
		\label{ineq:Ja2_grad}
	\end{equation}
	Therefore, from \eqref{prf:lem:terror_cont_norm:grad01}, \eqref{ineq:R_grad}, and \eqref{ineq:Ja2_grad}, we obtain \eqref{lem:terror_cont_norm:grad}. 
	
	Finally, we prove \eqref{lem:terror_cont_norm:lap}.	
	From \eqref{lh_def} and \eqref{taylor01} with $\ell=\WeightOrder+3$, we have 
	\begin{align}
		\LapApp{} \f(\x) 
		& = 2\Dim \sum_{1\leq |\multiindexa|\leq\WeightOrder+2} \frac{D^{\multiindexa}\f(\x)}{\multiindexa!} \sum_{\i\in\IndexSetXHAst{\x}{\h}} \Pv{\i} \frac{(\x-\Pt{\i})^{\multiindexa}}{|\x-\Pt{\i}|^2} \wh(|\x-\Pt{\i}|)
		\nn\\
		& \quad + 2\Dim \sum_{|\multiindexa|=\WeightOrder+3}\sum_{\i\in\IndexSetXHAst{\x}{\h}}R_{\multiindexa}(\Pt{\i},\x)\Pv{\i} \frac{(\x-\Pt{\i})^{\multiindexa}}{|\x-\Pt{\i}|^2} \wh(|\x-\Pt{\i}|). 
	\end{align}
	From \eqref{eq:d_int_y^b_y^-2dy}, we have
	\begin{align}
		2\Dim\sum_{|\multiindexa|=2}\frac{D^{\multiindexa}\f(\x)}{\multiindexa!} \int_{\dRd}\frac{\y^\multiindexa}{|\y|^2} \wh(|\y|)\dy = \Lap\f(\x). 
	\end{align}
	Hypothesis \ref{hypo:ref_weight:order} with $\WeightOrder$ yields
	\begin{align}
		\int_{\dRd} \frac{\y^\multiindexa}{|\y|^2} \wh(|\y|)\dy = 0,\qquad\mbox{$\multiindexa\in\dMI{\Dim}$ with $|\multiindexa|=1$ or $ 3\leq|\multiindexa|\leq\WeightOrder+2$}. 
	\end{align}
	Therefore, we have
	\begin{align}
		\Lap\f(\x) - \LapApp{}\f(\x)
		&=2\Dim\sum_{1\leq|\multiindexa|\leq\WeightOrder+2} \frac{D^{\multiindexa}\f(\x)}{\multiindexa!} \errorfuncB{\multiindexa}{2}(\x)
		\\
		&\quad+ 2\Dim\sum_{|\multiindexa|=\WeightOrder+3}\sum_{\i\in\IndexSetXHAst{\x}{\h}}R_{\multiindexa}(\Pt{\i},\x)\Pv{\i} \frac{(\x-\Pt{\i})^{\multiindexa}}{|\x-\Pt{\i}|^2}\wh(|\x-\Pt{\i}|). 
		\label{prf:lem:estimate:lap:01}
	\end{align}
	From \eqref{ineq:R_a[f](y_x)}, we have
	\begin{equation}
		\BrA{\sum_{|\multiindexa|=\WeightOrder+3}\sum_{\i\in\IndexSetXHAst{\x}{\h}}R_{\multiindexa}(\Pt{\i},\x)\Pv{\i} \frac{(\x-\Pt{\i})^{\multiindexa}}{|\x-\Pt{\i}|^2}\wh(|\x-\Pt{\i}|)}
		\leq\gconst|\errorfuncC{\WeightOrder+1}(\x)|\SNormC{\f}{\WeightOrder+3}{\DmHOl}. 
		\label{eq:int_y^a_wh(|y|)dy_lap}
	\end{equation}
	Moreover, we have 
	\begin{align}
		\BrA{\sum_{1\leq |\multiindexa|\leq\WeightOrder+2} \frac{D^{\multiindexa}\f(\x)}{\multiindexa!} \errorfuncB{\multiindexa}{2}(\x)}
		&\leq\gconst\NormC{\f}{\WeightOrder+2}{\DmOl}\sum_{1\leq|\multiindexa|\leq\WeightOrder+2}\BrA{\errorfuncB{\multiindexa}{2}(\x)}.
		\label{ineq:Ja2_lap}
	\end{align}
	Therefore, from \eqref{prf:lem:estimate:lap:01}, \eqref{eq:int_y^a_wh(|y|)dy_lap}, and \eqref{ineq:Ja2_lap}, we obtain \eqref{lem:terror_cont_norm:lap}. 
\end{proof}

Next, we show estimates of $\errorfuncA{\multiindexa}$, $\errorfuncB{\multiindexa}{\intA}$, and $\errorfuncC{\intA}$. 

\begin{lemma}
	\label{lem:Ja}
	There exists a positive constant $\gconst$ independent of $\N$ such that
	\begin{align}
		\BrNNd{\errorfuncA{\multiindexa}}_{\FsCz{\DmOl}} \leq \gconst \BrS{1+2\frac{\IndPtSet}{\h}}^{\Dim}\BrS{\frac{\IndPtSet+\IndPvSet}{\h}},\qquad \multiindexa\in\dMI{\Dim}. 
		\label{lem:ineq:Ja}
	\end{align}
\end{lemma}

\begin{proof}
	We arbitrarily fix $\x\in\DmOl$, $\multiindexa\in\dMI{\Dim}$, and particle volume decomposition $\PvDmSet=\{\PvDm{\i}\mid\i=1,2,\dots,\N\}$ and split $\errorfuncA{\multiindexa}$ into
	\begin{equation}
		\errorfuncA{\multiindexa}(\x) = \prffuncA{1}(\x) + \prffuncA{2}(\x) + \prffuncA{3}(\x)
	\end{equation}
	with
	\begin{align}
		\prffuncA{1}(\x) &\deq \sum_{\i\in\IndexSetXH{\x}{\h}}\Pv{\i} (\Pt{\i}-\x)^{\multiindexa}\wh(|\x-\Pt{\i}|)-\sum_{\i=1}^{\N}\sum_{\j=1}^{\N} \Meas{\VoroCell{\j}\cap\PvDm{\i}} (\Pt{\i}-\x)^{\multiindexa}\wh(|\x-\Pt{\j}|), 
		\\
		\prffuncA{2}(\x) &\deq \sum_{\i=1}^{\N}\sum_{\j=1}^{\N}(\Pt{\i}-\x)^{\multiindexa}\int_{\VoroCell{\j}\cap\PvDm{\i}}\{\wh(|\x-\Pt{\j}|)-\wh(|\x-\y|)\}\dy,
		\\
		\prffuncA{3}(\x) &\deq\sum_{\i=1}^{\N}\sum_{\j=1}^{\N}(\Pt{\i}-\x)^{\multiindexa}\int_{\VoroCell{\j}\cap\PvDm{\i}}\wh(|\x-\y|)\dy-\int_{\dRd}\y^\multiindexa \wh(|\y|)\dy. 
	\end{align}	
	Then, we estimate $\prffuncA{1}$, $\prffuncA{2}$, and $\prffuncA{3}$. 
	
	First, we estimate $\prffuncA{1}$. 
	Because 
	\begin{equation}
		\sum_{\j=1}^{\N} \Meas{\VoroCell{\j}\cap\PvDm{\i}}=\Pv{\i},\qquad \i=1,2,\dots,\N, 
		\label{sum_|VoroCell_PvDm|}
	\end{equation}
	we can rewrite $\prffuncA{1}$ as
	\begin{equation}
		\prffuncA{1}=\sum_{\i=1}^{\N}\sum_{\j=1}^{\N}\Meas{\VoroCell{\j}\cap\PvDm{\i}} (\Pt{\i}-\x)^{\multiindexa}\{\wh(|\x-\Pt{\i}|)-\wh(|\x-\Pt{\j}|)\}. 
	\end{equation}
	From
	\begin{equation}
		|(\y-\x)^{\multiindexa}| \leq \diam(\DmH)^{|\multiindexa|}, \qquad\y\in\DmH, 
		\label{y-x_leq_diam}
	\end{equation}
	we obtain
	\begin{align}
		\BrA{\prffuncA{1}(\x)} &\leq\gconst\sum_{\i=1}^{\N}\sum_{\j=1}^{\N}\Meas{\VoroCell{\j}\cap\PvDm{\i}}|\wh(|\x-\Pt{\i}|)-\wh(|\x-\Pt{\j}|)|. 
		\label{prf:estimate:E1:01}
	\end{align}
	From 
	\begin{equation}
		|\wh(|\x-\y|)-\wh(|\x-\z|)| = 0, \qquad \fa\y, \z\in\dRd\setminus\OBall{\x}{\h},
		\label{eq:|wi-wj|=0}
	\end{equation}
	we have 
	\begin{multline}
		\sum_{\i=1}^{\N}\sum_{\j=1}^{\N}\Meas{\VoroCell{\j}\cap\PvDm{\i}}|\wh(|\x-\Pt{\i}|)-\wh(|\x-\Pt{\j}|)|\\
		\leq\sum_{\i\in\IndexSetXH{\x}{\h}}\sum_{\j=1}^{\N}\Meas{\VoroCell{\j}\cap\PvDm{\i}}|\wh(|\x-\Pt{\i}|)-\wh(|\x-\Pt{\j}|)|
		+\sum_{\i=1}^{\N}\sum_{\j\in\IndexSetXH{\x}{\h}}\Meas{\VoroCell{\j}\cap\PvDm{\i}}|\wh(|\x-\Pt{\i}|)-\wh(|\x-\Pt{\j}|)|\\
		=\sum_{\i\in\IndexSetXH{\x}{\h}}\sum_{\j=1}^{\N}(\Meas{\VoroCell{\i}\cap\PvDm{\j}}+\Meas{\VoroCell{\j}\cap\PvDm{\i}})|\wh(|\x-\Pt{\i}|)-\wh(|\x-\Pt{\j}|)|.
		\label{prf:E1_estimate:lem:01}
	\end{multline}
	Because $\wh$ is absolutely continuous, we have
	\begin{align}
		\BrA{\wh(|\x-\y|)-\wh(|\x-\z|)} 
		&= \BrA{\{(\x-\y)-(\x-\z)\}\int_0^1\dfrac{d}{dr}\wh(t |\x-\y|+(1-t)|\x-\z|)dt}
		\nn\\
		&\leq |\y-\z| \BrA{  \int_0^1\dfrac{d}{dr}\wh(t |\x-\y|+(1-t)|\x-\z|)dt}
		\nn\\
		&\leq |\y-\z|  \int_0^h\BrA{\dfrac{d}{dr}\wh(r)}dr
		\nn\\
		&\leq \dfrac{|\y-\z|}{\h^{\Dim+1}}\int_0^1\BrA{\dfrac{d}{dr}\w(r)}dr,
		\label{ineq:|e_whx|}
	\end{align}
	for all $\y, \z\in\dRd$. 
	Moreover, we have
	\begin{align}
		\sum_{\i\in\IndexSetXH{\x}{r}}|\VoroCell{\i}|\leq\Meas{\OBall{\x}{1}}\BrS{r+\IndPtSet}^{\Dim},\qquad \fa r\in\dRpz. 
		\label{Sum_VoroCell_leq}
	\end{align}
	From \eqref{prf:E1_estimate:lem:01}, \eqref{ineq:|e_whx|}, and \eqref{Sum_VoroCell_leq}, we have
	\begin{align}
		\sum_{\i=1}^{\N}\sum_{\j=1}^{\N}\Meas{\VoroCell{\j}\cap\PvDm{\i}}|\wh(|\x-\Pt{\i}|)-\wh(|\x-\Pt{\j}|)|&\leq\dfrac{\gconst}{\h^{\Dim+1}}\sum_{\i\in\IndexSetXH{\x}{\h}}\sum_{\j=1}^{\N}(\Meas{\VoroCell{\i}\cap\PvDm{\j}}+\Meas{\VoroCell{\j}\cap\PvDm{\i}})|\Pt{\i}-\Pt{\j}|
		\nn\\
		&\leq\dfrac{\gconst}{\h^{\Dim+1}}\sum_{\i\in\IndexSetXH{\x}{\h}} |\VoroCell{\i}|\sum_{\j=1}^{\N} \frac{\Meas{\VoroCell{\i}\cap\PvDm{\j}}+\Meas{\VoroCell{\j}\cap\PvDm{\i}}}{|\VoroCell{\i}|}|\Pt{\i}-\Pt{\j}|
		\nn\\
		&\leq\gconst\dfrac{\IndPvSetLocal{\PvDmSet}}{\h^{\Dim+1}}\sum_{\i\in\IndexSetXH{\x}{\h}}|\VoroCell{\i}|
		\nn\\
		&\leq\gconst\BrS{1+\frac{\IndPtSet}{\h}}^{\Dim} \frac{\IndPvSetLocal{\PvDmSet}}{\h}.
		\label{prf:estimate:E1:02}
	\end{align}
	Therefore, from \eqref{prf:estimate:E1:01} and \eqref{prf:estimate:E1:02}, we obtain
	\begin{align}
		\BrA{\prffuncA{1}(\x)}
		& \leq\gconst\BrS{1+\frac{\IndPtSet}{\h}}^{\Dim} \frac{\IndPvSetLocal{\PvDmSet}}{\h}. 
		\label{prf:E1_estimate}
	\end{align}
	
	Next, we estimate $\prffuncA{2}$. 
	Because $\supp(\wh) = [0,\h]$ and $\VoroCell{\j}\subset\OBall{\Pt{\j}}{\IndPtSet}$, we have
	\begin{equation}
		\int_{\VoroCell{\j}\cap\PvDm{\i}}|\wh(|\x-\Pt{\j}|)-\wh(|\x-\y|)|\dy=0,\qquad\i=1,2,\dots,\N,~\j\not\in\IndexSetXH{\x}{\h+\IndPtSet}. 
		\label{prf:E2_estimate:01}
	\end{equation}
	From \eqref{prf:E2_estimate:01}, we have
	\renewcommand{\Space}{\quad}
	\begin{align}
		\sum_{\i=1}^{\N}\sum_{\j=1}^{\N}\int_{\VoroCell{\j}\cap\PvDm{\i}}|\wh(|\x-\Pt{\j}|)-\wh(|\x-\y|)|\dy&=\sum_{\i=1}^{\N}\sum_{\j\in\IndexSetXH{\x}{\h+\IndPtSet}}\int_{\VoroCell{\j}\cap\PvDm{\i}}|\wh(|\x-\Pt{\j}|)-\wh(|\x-\y|)|\dy\nn\\
		&=\sum_{\j\in\IndexSetXH{\x}{\h+\IndPtSet}}\int_{\VoroCell{\j}}|\wh(|\x-\Pt{\j}|)-\wh(|\x-\y|)|\dy.
	\end{align}
	Moreover, from \eqref{ineq:|e_whx|} and \eqref{Sum_VoroCell_leq}, we have
	\begin{align}
		\sum_{\i=1}^{\N}\sum_{\j=1}^{\N}\int_{\VoroCell{\j}\cap\PvDm{\i}}|\wh(|\x-\Pt{\j}|)-\wh(|\x-\y|)|\dy
		&\leq\dfrac{\gconst}{\h^{\Dim+1}}\sum_{\j\in\IndexSetXH{\x}{\h+\IndPtSet}}\int_{\VoroCell{\j}}|\Pt{\j}-\y|\dy
		\nn\\
		&\leq\dfrac{\gconst}{\h^{\Dim+1}}\sum_{\j\in\IndexSetXH{\x}{\h+\IndPtSet}}\int_{\VoroCell{\j}}|\Pt{\j}-\y|\dy
		\nn\\
		&\leq\gconst\dfrac{\IndPtSet}{\h^{\Dim+1}}\sum_{\j\in\IndexSetXH{\x}{\h+\IndPtSet}}|\VoroCell{\j}|
		\nn\\
		&\leq\gconst\BrS{1+2\dfrac{\IndPtSet}{\h}}^{\Dim}\dfrac{\IndPtSet}{\h}. 
		\label{prf:estimate:E2:01}
	\end{align}
	Therefore, from \eqref{y-x_leq_diam} and \eqref{prf:estimate:E2:01}, we obtain
	\begin{align}
		|\prffuncA{2}(\x)|
		&\leq\sum_{\i=1}^{\N}\sum_{\j=1}^{\N}\BrA{(\Pt{\i}-\x)^{\multiindexa}}\int_{\VoroCell{\j}\cap\PvDm{\i}}\BrA{\wh(|\x-\Pt{\j}|)-\wh(|\x-\y|)}\dy\nn\\
		&\leq\gconst \sum_{\i=1}^{\N}\sum_{\j=1}^{\N}\int_{\VoroCell{\j}\cap\PvDm{\i}}\BrA{\wh(|\x-\Pt{\j}|)-\wh(|\x-\y|)}\dy\nn\\
		&\leq\gconst\BrS{1+2\dfrac{\IndPtSet}{\h}}^{\Dim}\dfrac{\IndPtSet}{\h}. 
	\end{align}
	
	Finally, we estimate $\prffuncA{3}$. 
	Because
	\begin{equation}
		\int_{\dRd}\y^\multiindexa \wh(|\y|)\dy=\int_{\DmH}(\y-\x)^\multiindexa \wh(|\x-\y|)\dy,
	\end{equation}
	we can rewrite $\prffuncA{3}$ as
	\begin{align}
		\prffuncA{3}(\x) &=\sum_{\i=1}^{\N}\sum_{\j=1}^{\N}\int_{\VoroCell{\j}\cap\PvDm{\i}}\{(\Pt{\i}-\x)^{\multiindexa}-(\y-\x)^\multiindexa\}\wh(|\x-\y|)\dy. 
	\end{align}
	Because $\prffuncA{3}=0$ when $|\multiindexa|=0$, we estimate when $|\multiindexa|\geq1$. 
	Let $\multiindexb_k\,(k=1,2,\dots,|\multiindexa|)$ be $\Dim$-dimensional multi-indices with satisfying 
	\begin{equation}
		\sum_{k=1}^{|\multiindexa|}\multiindexb_k=\multiindexa,\qquad|\multiindexb_k|=1~(k=1,2,\dots,|\multiindexa|). 
	\end{equation}
	Then, we have, for all $y,z\in\dRd$,
	\begin{align}
		\BrA{\y^{\multiindexa}-\z^{\multiindexa}} 
		& \leq \BrA{\y^{\multiindexa}-\y^{\multiindexa-\multiindexb_1}\z^{\multiindexb_1}} + \BrA{\y^{\multiindexa-\multiindexb_1}\z^{\multiindexb_1}-\z^{\multiindexa}}
		\nn\\
		& \leq \BrA{\y-\z} |\y|^{|\multiindexa|-1}+ \BrA{\y^{\multiindexa-\multiindexb_1}-\z^{\multiindexa-\multiindexb_1}}|\z|
		\nn\\
		& \leq \BrA{\y-\z} |\y|^{|\multiindexa|-1}+\BrA{\y-\z} |\y|^{|\multiindexa|-2}|\z|+ \BrA{\y^{\multiindexa-\multiindexb_1-\multiindexb_2}-\z^{\multiindexa-\multiindexb_1-\multiindexb_2}}|\z|^2
		\nn\\
		&~\vdots
		\nn\\
		&\leq\BrA{\y-\z}\sum_{k=1}^{|\multiindexa|}|\y|^{|\multiindexa|-k}|\z|^{k-1}. 
		\label{ineq:|y^alpha-z^alpha|}
	\end{align}
	From \eqref{y-x_leq_diam} and \eqref{ineq:|y^alpha-z^alpha|}, we obtain 
	\begin{align}
		|\prffuncA{3}(\x)|
		&\leq\sum_{\i=1}^{\N}\sum_{\j=1}^{\N}\int_{\VoroCell{\j}\cap\PvDm{\i}}|(\Pt{\i}-\x)^{\multiindexa}-(\y-\x)^\multiindexa||\wh(|\x-\y|)|\dy\nn\\
		&\leq\gconst\sum_{\i=1}^{\N}\sum_{\j=1}^{\N}\int_{\VoroCell{\j}\cap\PvDm{\i}}|\y-\Pt{\i}||\wh(|\x-\y|)|\dy.
		\label{prf:estimate:E3:01}
	\end{align}
	By $\supp(\wh) = [0,\h]$ and $\VoroCell{\j}\subset\OBall{\Pt{\j}}{\IndPtSet}$, if $\j\not\in\IndexSetXH{\x}{\h+\IndPtSet}$, then 
	\begin{equation}
		\int_{\VoroCell{\j}\cap\PvDm{\i}}|\y-\Pt{\i}||\wh(|\x-\y|)|\dy=0,\qquad\i=1,2,\dots,\N. 
		\label{prf:estimate:E3:02}
	\end{equation}
	Moreover, from $\w\in\FsWeightFunc\subset\FsCz{\dRpz}$, we have
	\begin{equation}
		|\wh(|\x-\y|)|=\dfrac{1}{\h^\Dim}\BrA{\w\BrS{\dfrac{|\x-\y|}{\h}}}\leq\dfrac{1}{\h^\Dim}\NormCz{\w}{\dRpz},\qquad \fa \y\in\DmH. 
		\label{prf:estimate:E3:03}
	\end{equation}
	From \eqref{Sum_VoroCell_leq}, \eqref{prf:estimate:E3:02}, and \eqref{prf:estimate:E3:03}, we have
	\begin{align}
		&\sum_{\i=1}^{\N}\sum_{\j=1}^{\N}\int_{\VoroCell{\j}\cap\PvDm{\i}}|\y-\Pt{\i}||\wh(|\x-\y|)|\dy
		\nn\\
		&\quad=\sum_{\i=1}^{\N}\sum_{\j\in\IndexSetXH{\x}{\h+\IndPtSet}}\int_{\VoroCell{\j}\cap\PvDm{\i}}|\y-\Pt{\i}||\wh(|\x-\y|)|\dy
		\nn\\
		&\quad\leq\dfrac{\gconst}{\h^\Dim}\sum_{\i=1}^{\N}\sum_{\j\in\IndexSetXH{\x}{\h+\IndPtSet}}\int_{\VoroCell{\j}\cap\PvDm{\i}}|\y-\Pt{\i}|\dy
		\nn\\
		&\quad\leq\dfrac{\gconst}{\h^\Dim}\sum_{\i=1}^{\N}\sum_{\j\in\IndexSetXH{\x}{\h+\IndPtSet}}\int_{\VoroCell{\j}\cap\PvDm{\i}}(|\y-\Pt{\j}|+|\Pt{\j}-\Pt{\i}|)\dy\nn\\
		&\quad\leq\dfrac{\gconst}{\h^\Dim}\BrS{\IndPtSet\sum_{\j\in\IndexSetXH{\x}{\h+\IndPtSet}}\Meas{\VoroCell{\j}}+\sum_{\j\in\IndexSetXH{\x}{\h+\IndPtSet}}\sum_{\i=1}^{\N}\Meas{\VoroCell{\j}\cap\PvDm{\i}}|\Pt{\j}-\Pt{\i}|}\nn\\
		&\quad\leq\dfrac{\gconst}{\h^\Dim}\BrS{\sum_{\j\in\IndexSetXH{\x}{\h+\IndPtSet}}\Meas{\VoroCell{\j}}}\BrM{\IndPtSet+\max_{\j=1,2,\dots,\N}\BrS{\sum_{\i=1}^{\N}\dfrac{\Meas{\VoroCell{\i}\cap\PvDm{\j}}+\Meas{\VoroCell{\j}\cap\PvDm{\i}}}{\Meas{\VoroCell{\j}}}|\Pt{\j}-\Pt{\i}|}}\nn\\
		&\quad\leq\gconst\BrS{1+\dfrac{\IndPtSet}{\h}}^{\Dim}\BrS{\IndPtSet+\IndPvSetLocal{\PvDmSet}}.
		\label{prf:estimate:E3:04}
	\end{align}
	Therefore, from \eqref{prf:estimate:E3:01}, \eqref{prf:estimate:E3:04}, and $\h\leq\H$, we obtain
	\begin{align}
		|\prffuncA{3}(\x)| &\leq\gconst\BrS{1+2\dfrac{\IndPtSet}{\h}}^{\Dim}(\IndPtSet+\IndPvSetLocal{\PvDmSet})\nn\\
		&\leq\gconst\BrS{1+2\dfrac{\IndPtSet}{\h}}^{\Dim}\frac{\IndPtSet+\IndPvSetLocal{\PvDmSet}}{\h}.
	\end{align}
	From the estimates of $\prffuncA{1}$, $\prffuncA{2}$, and $\prffuncA{3}$, we obtain
	\begin{equation}
		\BrNNd{\errorfuncA{\multiindexa}}_{\FsCz{\DmOl}} \leq \gconst \BrS{1+2\frac{\IndPtSet}{\h}}^{\Dim}\frac{\IndPtSet+\IndPvSetLocal{\PvDmSet}}{\h}. 
	\end{equation}
	Because $\PvDmSet$ is arbitrary, we establish \eqref{lem:ineq:Ja}. 
\end{proof}

\begin{lemma}
	\label{lem:Jam}
	Suppose that reference weight function $\w$ satisfies Hypothesis \ref{hypo:ref_weight:smooth} with $\WeightOrderSmooth$. 
	Then, there exists a positive constant $\gconst$ independent of $\N$ such that for all $\multiindexa\in\dMI{\Dim}$ and $\intA\in\dN$ with $1\leq\intA-\WeightOrderSmooth\leq|\multiindexa|$, 
	\begin{align}
		\BrNNd{\errorfuncB{\multiindexa}{\intA}}_{\FsCz{\DmOl}} \leq \gconst \BrS{1+2\frac{\IndPtSet}{\h}}^{\Dim} \frac{\IndPtSet+\IndPvSet}{\h^{\WeightOrderSmooth+1}}. 
		\label{lem:ineq:Jam}
	\end{align}
\end{lemma}

\begin{proof}
	We arbitrarily fix $\x\in\DmOl$, $\multiindexa\in\dMI{\Dim}$, particle volume decomposition $\PvDmSet=\{\PvDm{\i}\mid\i=1,2,\dots,\N\}$, and $\intA\in\dN$ with $1\leq\intA-\WeightOrderSmooth\leq|\multiindexa|$ and
	split $\errorfuncB{\multiindexa}{\intA}$ into
	\begin{equation}
		\errorfuncB{\multiindexa}{\intA}(\x) = \prffuncA{4}(\x) + \prffuncA{5}(\x) + \prffuncA{6}(\x)
	\end{equation}
	with
	\renewcommand{\Space}{\quad}
	\begin{align}
		\prffuncA{4}(\x) &\deq \sum_{\i\in\IndexSetXHAst{\x}{\h}}\Pv{\i}\frac{(\Pt{\i}-\x)^{\multiindexa}}{|\Pt{\i}-\x|^{\intA}} \wh(|\x-\Pt{\i}|)
		-\sum_{\i\in\IndexSetXHAst{\x}{\infty}}\sum_{\j\in\IndexSetXHAst{\x}{\infty}}\Meas{\VoroCell{\j}\cap\PvDm{\i}}\frac{(\Pt{\i}-\x)^{\multiindexa}}{|\Pt{\i}-\x|^{\intA-\WeightOrderSmooth}}\dfrac{\wh(|\x-\Pt{\j}|)}{|\x-\Pt{\j}|^{\WeightOrderSmooth}},
		\\
		\prffuncA{5}(\x)&\deq\sum_{\i\in\IndexSetXHAst{\x}{\infty}}\sum_{\j\in\IndexSetXHAst{\x}{\infty}}\Meas{\VoroCell{\j}\cap\PvDm{\i}}\frac{(\Pt{\i}-\x)^{\multiindexa}}{|\Pt{\i}-\x|^{\intA-\WeightOrderSmooth}}\dfrac{\wh(|\x-\Pt{\j}|)}{|\x-\Pt{\j}|^{\WeightOrderSmooth}}-\sum_{\i\in\IndexSetXHAst{\x}{\infty}}\sum_{\j=1}^{\N}\frac{(\Pt{\i}-\x)^{\multiindexa}}{|\Pt{\i}-\x|^{\intA-\WeightOrderSmooth}}\int_{\VoroCell{\j}\cap\PvDm{\i}}\dfrac{\wh(|\x-\y|)}{|\x-\y|^{\WeightOrderSmooth}}\dy,
		\\
		\prffuncA{6}(\x) &\deq\sum_{\i\in\IndexSetXHAst{\x}{\infty}}\sum_{\j=1}^{\N}\frac{(\Pt{\i}-\x)^{\multiindexa}}{|\Pt{\i}-\x|^{\intA-\WeightOrderSmooth}}\int_{\VoroCell{\j}\cap\PvDm{\i}}\dfrac{\wh(|\x-\y|)}{|\x-\y|^{\WeightOrderSmooth}}\dy-\int_{\dRd} \frac{\y^\multiindexa}{|\y|^\intA} \wh(|\y|)\dy. 
	\end{align}	
	Then, we estimate $\prffuncA{4}$, $\prffuncA{5}$, and $\prffuncA{6}$. 
	
	First, we estimate $\prffuncA{4}$ and set $\wn{k}$ as \eqref{def:wn} and $\wnh{k}$ as
	\begin{align}
		\wnh{k}(r) & \deq \frac{1}{\h^{\Dim+k}}\wn{k}\left(\frac{r}{\h}\right), \qquad r\in\dRpz.
		\label{def:wnh}
	\end{align}
	Then, from \eqref{sum_|VoroCell_PvDm|}, we can rewrite $\prffuncA{4}$ as
	\begin{equation}
		\prffuncA{4}(\x) =
		\sum_{\i\in\IndexSetXHAst{\x}{\infty}}\sum_{\j=1}^{\N}\Meas{\VoroCell{\j}\cap\PvDm{\i}}\frac{(\Pt{\i}-\x)^{\multiindexa}}{|\Pt{\i}-\x|^{\intA-\WeightOrderSmooth}}\{\wnh{\WeightOrderSmooth}(|\x-\Pt{\i}|)-\wnh{\WeightOrderSmooth}(|\x-\Pt{\j}|)\}.
	\end{equation}
	Because 
	\begin{align}
		\BrA{\frac{(\Pt{\i}-\x)^{\multiindexa}}{|\Pt{\i}-\x|^{\intA-\WeightOrderSmooth}}}\leq |\Pt{\i}-\x|^{|\multiindexa|-\intA+\WeightOrderSmooth}\leq \diam(\DmH)^{|\multiindexa|-\intA+\WeightOrderSmooth}, \quad \i\in\IndexSetXHAst{\x}{\infty}, 
		\label{ineq:|psi_yz|}
	\end{align}
	we obtain 
	\begin{equation}
		|\prffuncA{4}(\x)|\leq\gconst\sum_{\i=1}^{\N}\sum_{\j=1}^{\N}\Meas{\VoroCell{\j}\cap\PvDm{\i}}\BrA{\wnh{\WeightOrderSmooth}(|\x-\Pt{\i}|)-\wnh{\WeightOrderSmooth}(|\x-\Pt{\j}|)}. 
	\end{equation}
	From $\supp(\wnh{\WeightOrderSmooth}) = [0,\h]$, we have
	\begin{align}
		\wnh{\WeightOrderSmooth}(|\x-\Pt{\i}|)-\wnh{\WeightOrderSmooth}(|\x-\Pt{\j}|)= 0,\qquad \i, \j\not\in\IndexSetXHAst{\x}{\h}. 
		\label{eq:e_wh^n=0}
	\end{align}
	Thus, we obtain
	\begin{multline}
		|\prffuncA{4}(\x)|
		\leq\gconst\Bigg(\sum_{\i\in\IndexSetXHAst{\x}{\h}}\sum_{\j=1}^{\N}\Meas{\VoroCell{\j}\cap\PvDm{\i}}\BrA{\wnh{\WeightOrderSmooth}(|\x-\Pt{\i}|)-\wnh{\WeightOrderSmooth}(|\x-\Pt{\j}|)}\\
		\qquad+\sum_{\i=1}^{\N}\sum_{\j\in\IndexSetXHAst{\x}{\h}}\Meas{\VoroCell{\j}\cap\PvDm{\i}}\BrA{\wnh{\WeightOrderSmooth}(|\x-\Pt{\i}|)-\wnh{\WeightOrderSmooth}(|\x-\Pt{\j}|)}\Bigg).
		\label{prf:estimate:E4:01}
	\end{multline}
	Using an argument similar to \eqref{ineq:|e_whx|}, if $\w$ satisfies Hypothesis \ref{hypo:ref_weight:smooth} with $\WeightOrderSmooth$, then for all $\y, \z\in\dRd$, 
	\begin{align}
		\BrANd{\wnh{\WeightOrderSmooth}(|\x-\y|)-\wnh{\WeightOrderSmooth}(|\x-\z|)} 
		& \leq  \dfrac{|\y-\z|}{\h^{\Dim+k+1}}\int_0^1\BrA{ \dfrac{d}{dr}\wn{\WeightOrderSmooth}(r)}dr. 
		\label{ineq:|e_whx^n|}
	\end{align}
	From \eqref{prf:estimate:E4:01} and \eqref{ineq:|e_whx^n|}, we obtain
	\begin{align}
		|\prffuncA{4}(\x)|
		&\leq\dfrac{\gconst}{\h^{\Dim+\WeightOrderSmooth+1}}\sum_{\i\in\IndexSetXHAst{\x}{\h}}\sum_{\j=1}^{\N}\BrS{\Meas{\VoroCell{\i}\cap\PvDm{\j}}+\Meas{\VoroCell{\j}\cap\PvDm{\i}}}\BrA{\Pt{\i}-\Pt{\j}}
		\nn\\
		&\leq\dfrac{\gconst}{\h^{\Dim+\WeightOrderSmooth+1}}\sum_{\i\in\IndexSetXHAst{\x}{\h}}\Meas{\VoroCell{\i}}\sum_{\j=1}^{\N}\dfrac{\Meas{\VoroCell{\i}\cap\PvDm{\j}}+\Meas{\VoroCell{\j}\cap\PvDm{\i}}}{\Meas{\VoroCell{\i}}}\BrA{\Pt{\i}-\Pt{\j}}
		\nn\\
		&\leq\gconst\BrS{1+\frac{\IndPtSet}{\h}}^{\Dim} \frac{\IndPvSetLocal{\PvDmSet}}{\h^{\WeightOrderSmooth+1}}. 
	\end{align}
	
	Next, we estimate $\prffuncA{5}$. 
	By using $\wnh{\WeightOrderSmooth}$, we can rewrite $\prffuncA{5}$ as
	\begin{equation}
		\prffuncA{5}(\x)=\sum_{\i\in\IndexSetXHAst{\x}{\infty}}\sum_{\j=1}^\N\frac{(\Pt{\i}-\x)^{\multiindexa}}{|\Pt{\i}-\x|^{\intA-\WeightOrderSmooth}}\int_{\VoroCell{\j}\cap\PvDm{\i}}\BrM{\wnh{\WeightOrderSmooth}(|\x-\Pt{\j}|)-\wnh{\WeightOrderSmooth}(|\x-\y|)}\dy. 
	\end{equation}
	From \eqref{ineq:|psi_yz|}, we obtain 
	\begin{align}
		|\prffuncA{5}(\x)|
		&\leq\gconst\sum_{\i=1}^\N\sum_{\j=1}^\N\int_{\VoroCell{\j}\cap\PvDm{\i}}\BrA{\wnh{\WeightOrderSmooth}(|\x-\Pt{\j}|)-\wnh{\WeightOrderSmooth}(|\x-\y|)}\dy\nn\\
		&\leq\gconst\sum_{\j=1}^\N\int_{\VoroCell{\j}}\BrA{\wnh{\WeightOrderSmooth}(|\x-\Pt{\j}|)-\wnh{\WeightOrderSmooth}(|\x-\y|)}\dy. 
	\end{align}
	By $\supp(\wnh{\WeightOrderSmooth}) = [0,\h]$ and $\VoroCell{\j}\subset\OBall{\Pt{\j}}{\IndPtSet}$, we have
	\begin{align}
		\int_{\VoroCell{\j}}\BrA{\wnh{\WeightOrderSmooth}(|\x-\Pt{\j}|)-\wnh{\WeightOrderSmooth}(|\x-\y|)}\dy= 0,\qquad \j\not\in\IndexSetXHAst{\x}{\h+\IndPtSet}. 
		\label{prf:estimate:E5:01}
	\end{align}
	From \eqref{ineq:|e_whx^n|} and \eqref{prf:estimate:E5:01}, we obtain
	\begin{align}
		|\prffuncA{5}(\x)|
		&\leq\gconst\sum_{\j\in\IndexSetXH{\x}{\h+\IndPtSet}}\int_{\VoroCell{\j}}\BrA{\wnh{\WeightOrderSmooth}(|\x-\Pt{\j}|)-\wnh{\WeightOrderSmooth}(|\x-\y|)}\dy
		\nn\\
		&\leq\dfrac{\gconst}{\h^{\Dim+\WeightOrderSmooth+1}}\int_{\VoroCell{\j}}\BrA{\Pt{\j}-\y}\dy
		\nn\\
		&\leq\gconst\dfrac{\IndPtSet}{\h^{\Dim+\WeightOrderSmooth+1}}\sum_{\j\in\IndexSetXH{\x}{\h+\IndPtSet}}\Meas{\VoroCell{\j}}
		\nn\\
		&\leq\gconst\BrS{1+2\frac{\IndPtSet}{\h}}^{\Dim} \dfrac{\IndPtSet}{\h^{\WeightOrderSmooth+1}}.
	\end{align}
	
	Finally, we estimate $\prffuncA{6}$. 
	Using $\wnh{\WeightOrderSmooth}$, we can rewrite $\prffuncA{6}$ as
	\begin{align}
		\prffuncA{6}(\x)&=\sum_{\i\in\IndexSetXH{\x}{\infty}}\sum_{\j=1}^{\N}\int_{\VoroCell{\j}\cap\PvDm{\i}}\BrM{\frac{(\Pt{\i}-\x)^{\multiindexa}}{|\Pt{\i}-\x|^{\intA-\WeightOrderSmooth}}-\frac{(\y-\x)^{\multiindexa}}{|\y-\x|^{\intA-\WeightOrderSmooth}}}\wnh{\WeightOrderSmooth}(|\x-\y|)\dy\nn\\
		&\quad-\sum_{\i=1}^{\N}\sum_{\j=1}^{\N}\int_{\VoroCell{\j}\cap\PvDmFunc{\i}{\x}}\frac{(\y-\x)^{\multiindexa}}{|\y-\x|^{\intA-\WeightOrderSmooth}}\wnh{\WeightOrderSmooth}(|\x-\y|)\dy, 
		\label{prf:estimate:E6:02}
	\end{align}
	where $\PvDmFunc{\i}{\x}$ is 
	\begin{equation}
		\PvDmFunc{\i}{\x} = 
		\begin{cases}
			\PvDm{\i},\qquad &\x=\Pt{\i},
			\\
			0,\qquad &\mbox{otherwize}. 
		\end{cases}
	\end{equation}
	For $\multiindexa\in\dMI{\Dim}$, let  $\multiindexb_j\,(j=1,2,\dots,|\multiindexa|)$ be $\Dim$-dimensional multi-indices satisfying
	\begin{align}
		|\multiindexb_j| = 1 \quad \mbox{and}  \quad \sum_{j=1}^{|\multiindexa|}\multiindexb_j = \multiindexa. 
	\end{align}
	Let $\multiindexb_j^\ast\,(j=0,1,\dots,|\multiindexa|)$ be $\Dim$-dimensional multi-indices defined as
	\begin{align}
		\multiindexb_j^\ast \deq
		\begin{cases}
			0,\qquad & j=0,
			\\
			\ds \sum_{\ell=1}^{j}\multiindexb_\ell,\qquad & j=1,2,\dots,|\multiindexa|. 
		\end{cases} 
	\end{align}
	For all $\y, \z \in\dRd\setminus\{0\}$, when $|\multiindexa|=\intA-\WeightOrderSmooth$, we have
	\begin{align}
		\BrA{\dfrac{\y^{\multiindexa}}{|\y|^{\intA-\WeightOrderSmooth}}-\dfrac{\z^{\multiindexa}}{|\z|^{\intA-\WeightOrderSmooth}}} 
		&\leq \sum_{j=0}^{\intA-\WeightOrderSmooth-1} \BrA{\dfrac{\y^{\multiindexb_{|\multiindexa|-j}^\ast}\z^{\multiindexb_{j}^\ast}}{|\y|^{\intA-\WeightOrderSmooth-j}|\z|^{j}}-\dfrac{\y^{\multiindexb_{|\multiindexa|-j-1}^\ast}\z^{\multiindexb_{j+1}^\ast}}{|\y|^{\intA-\WeightOrderSmooth-j-1}|\z|^{j+1}}} 
		\nn\\
		&\leq \sum_{j=0}^{\intA-\WeightOrderSmooth-1} \BrA{\dfrac{\y^{\multiindexb_{|\multiindexa|-j}^\ast}\z^{\multiindexb_{j}^\ast}-\y^{\multiindexb_{|\multiindexa|-j-1}^\ast}\z^{\multiindexb_{j+1}^\ast}}{|\y|^{\intA-\WeightOrderSmooth-j}|\z|^{j}}}+\sum_{j=0}^{\intA-\WeightOrderSmooth-1} \BrA{\dfrac{\y^{\multiindexb_{|\multiindexa|-j-1}^\ast}\z^{\multiindexb_{j+1}^\ast}}{|\y|^{\intA-\WeightOrderSmooth-j}|\z|^{j}}-\dfrac{\y^{\multiindexb_{|\multiindexa|-j-1}^\ast}\z^{\multiindexb_{j+1}^\ast}}{|\y|^{\intA-\WeightOrderSmooth-j-1}|\z|^{j+1}}}
		\nn\\
		&\leq 2(\intA-\WeightOrderSmooth)\dfrac{|\y-\z|}{|\y|}. 
		\label{prf:estimate:E6:03}
	\end{align}
	Moreover, from \eqref{ineq:|y^alpha-z^alpha|} and \eqref{prf:estimate:E6:03}, when $|\multiindexa|>\intA-\WeightOrderSmooth$, we have
	\begin{align}
		\BrA{\dfrac{\y^{\multiindexa}}{|\y|^{\intA-\WeightOrderSmooth}}-\dfrac{\z^{\multiindexa}}{|\z|^{\intA-\WeightOrderSmooth}}} 
		&\leq\BrA{\dfrac{\y^{\multiindexa}}{|\y|^{\intA-\WeightOrderSmooth}}-\dfrac{\y^{\multiindexb_{|\multiindexa|-\intA+\WeightOrderSmooth}^\ast}\z^{\multiindexb_{\intA-\WeightOrderSmooth}^\ast}}{|\z|^{\intA-\WeightOrderSmooth}}}+\BrA{\dfrac{\y^{\multiindexb_{|\multiindexa|-\intA+\WeightOrderSmooth}^\ast}\z^{\multiindexb_{\intA-\WeightOrderSmooth}^\ast}}{|\z|^{\intA-\WeightOrderSmooth}}-\dfrac{\z^{\multiindexa}}{|\z|^{\intA-\WeightOrderSmooth}}} 
		\nn\\
		&\leq|\y|^{|\multiindexa|-\intA+\WeightOrderSmooth}\BrA{\dfrac{\y^{\multiindexb_{\intA-\WeightOrderSmooth}^\ast}}{|\y|^{\intA-\WeightOrderSmooth}}-\dfrac{\z^{\multiindexb_{\intA-\WeightOrderSmooth}^\ast}}{|\z|^{\intA-\WeightOrderSmooth}}} +\BrA{\y^{\multiindexb_{|\multiindexa|-\intA+\WeightOrderSmooth}^\ast}-\z^{\multiindexb_{|\multiindexa|-\intA+\WeightOrderSmooth}^\ast}} 
		\nn\\
		&\leq 2(\intA-\WeightOrderSmooth)|\y-\z||\y|^{|\multiindexa|-\intA+\WeightOrderSmooth-1}+|\y-\z|\sum_{j=0}^{|\multiindexa|-\intA+\WeightOrderSmooth-1}|\y|^{j}|\z|^{|\multiindexa|-\intA+\WeightOrderSmooth-1-j}. 
	\end{align}
	Therefore, when $|\multiindexa|\geq\intA-\WeightOrderSmooth$, we have for all $\y\in\DmH\setminus\{\x\}$ and $\i\in\IndexSetXHAst{\x}{\infty}$, 
	\begin{align}
		\BrA{\frac{(\Pt{\i}-\x)^{\multiindexa}}{|\Pt{\i}-\x|^{\intA-\WeightOrderSmooth}}-\frac{(\y-\x)^{\multiindexa}}{|\y-\x|^{\intA-\WeightOrderSmooth}}}
		&\leq \gconst\dfrac{|\y-\Pt{\i}|}{|\y-\x|}. 
		\label{ineq:|psi(xy)-psy(xz)|_k=}
	\end{align}
	From \eqref{ineq:|psi(xy)-psy(xz)|_k=}, we obtain
	\begin{align}
		|\prffuncA{6}(\x)| &\leq\sum_{\i\neq\k}\sum_{\j=1}^{\N}\int_{\VoroCell{\j}\cap\PvDm{\i}}\BrA{\frac{(\Pt{\i}-\x)^{\multiindexa}}{|\Pt{\i}-\x|^{\intA-\WeightOrderSmooth}}-\frac{(\y-\x)^{\multiindexa}}{|\y-\x|^{\intA-\WeightOrderSmooth}}}\BrA{\wnh{\WeightOrderSmooth}(|\x-\y|)}\dy\nn\\
		&\quad+\BrA{\sum_{\i=1}^{\N}\sum_{\j=1}^{\N}\int_{\VoroCell{\j}\cap\PvDmFunc{\i}{\x}}\frac{(\y-\x)^{\multiindexa}}{|\y-\x|^{\intA-\WeightOrderSmooth}}\wnh{\WeightOrderSmooth}(|\x-\y|)\dy}\nn\\
		&\leq\gconst\sum_{\i\neq\k}\sum_{\j=1}^{\N} \int_{\VoroCell{\j}\cap\PvDm{\i}}|\y-\Pt{\i}|\BrA{\wnh{\WeightOrderSmooth+1}(|\x-\y|)}\dy\nn\\
		&\quad+\BrA{\sum_{\i=1}^{\N}\sum_{\j=1}^{\N}\int_{\VoroCell{\j}\cap\PvDmFunc{\i}{\x}}\frac{(\y-\x)^{\multiindexa}}{|\y-\x|^{\intA-\WeightOrderSmooth}}\wnh{\WeightOrderSmooth}(|\x-\y|)\dy}. 
		\label{prf:estimate:E6:05}
	\end{align}
	Because $|\multiindexa|\geq\intA-\WeightOrderSmooth$, we have
	\begin{align}
		\BrA{\sum_{\i=1}^{\N}\sum_{\j=1}^{\N}\int_{\VoroCell{\j}\cap\PvDmFunc{\i}{\x}}\frac{(\y-\x)^{\multiindexa}}{|\y-\x|^{\intA-\WeightOrderSmooth}}\wnh{\WeightOrderSmooth}(|\x-\y|)\dy}
		&\leq\gconst\sum_{\i=1}^{\N}\sum_{\j=1}^{\N}\int_{\VoroCell{\j}\cap\PvDmFunc{\i}{\x}}\BrA{\frac{(\y-\x)^{\multiindexa}}{|\y-\x|^{\intA-\WeightOrderSmooth-1}}}\BrA{\wnh{\WeightOrderSmooth+1}(|\x-\y|)}\dy
		\nn\\
		&\leq\gconst\sum_{\i=1}^{\N}\sum_{\j=1}^{\N}\int_{\VoroCell{\j}\cap\PvDmFunc{\i}{\x}}\BrA{\y-\x}\BrA{\wnh{\WeightOrderSmooth+1}(|\x-\y|)}\dy. 
	\end{align}
	Therefore, we have
	\begin{equation}
		|\prffuncA{6}(\x)| 
		\leq\gconst\sum_{\i=1}^{\N}\sum_{\j=1}^{\N} \int_{\VoroCell{\j}\cap\PvDm{\i}}|\y-\Pt{\i}|\BrA{\wnh{\WeightOrderSmooth+1}(|\x-\y|)}\dy. 
	\end{equation}
	Because for all $\y\in\DmH$, 
	\begin{equation}
		\BrA{\wnh{\WeightOrderSmooth+1}(|\x-\y|)}=\dfrac{1}{\h^{\Dim+\WeightOrderSmooth+1}}\BrA{\wn{\WeightOrderSmooth+1}\BrS{\dfrac{|\x-\y|}{\h}}}\leq\dfrac{1}{\h^{\Dim+\WeightOrderSmooth+1}}\NormCzNd{\wn{\WeightOrderSmooth+1}}{\dRpz}, 
	\end{equation}
	by the same procedure as \eqref{prf:estimate:E3:04}, we have
	\begin{equation}
		\sum_{\i=1}^{\N}\sum_{\j=1}^{\N} \int_{\VoroCell{\j}\cap\PvDm{\i}}|\y-\Pt{\i}|\BrA{\wnh{\WeightOrderSmooth+1}(|\x-\y|)}\dy\leq\gconst\BrS{1+\dfrac{\IndPtSet}{\h}}^{\Dim}\dfrac{\IndPtSet+\IndPvSetLocal{\PvDmSet}}{\h^{\WeightOrderSmooth+1}}.
	\end{equation}
	Therefore, we obtain
	\begin{align}
		|\prffuncA{6}(\x)|\leq\gconst\BrS{1+\dfrac{\IndPtSet}{\h}}^{\Dim}\dfrac{\IndPtSet+\IndPvSetLocal{\PvDmSet}}{\h^{\WeightOrderSmooth+1}}. 
	\end{align}
	
	From the estimates of $\prffuncA{4}$, $\prffuncA{5}$, and $\prffuncA{6}$, we obtain
	\begin{equation}
		\BrNNd{\errorfuncB{\multiindexa}{\intA}}_{\FsCz{\DmOl}}
		\leq\gconst\BrS{1+2\frac{\IndPtSet}{\h}}^{\Dim}\frac{\IndPtSet+\IndPvSetLocal{\PvDmSet}}{\h^{\WeightOrderSmooth+1}}. 
	\end{equation}
	Because $\PvDmSet$ is arbitrary, we establish \eqref{lem:ineq:Jam}. 
\end{proof}

\begin{lemma}
	\label{lem:Km}
	There exists a positive constant $\gconst$ independent of $\N$ such that  
	\begin{align}
		\NormCz{\errorfuncC{\intA}}{\DmOl} \leq \gconst\BrM{\BrS{1+2\frac{\IndPtSet}{\h}}^\Dim \dfrac{\IndPtSet+\IndPvSet}{\h}+\h^{\intA}},\quad \intA\in\dN. 
		\label{eq:lem:Km}
	\end{align}	
\end{lemma}

\begin{proof}
	We arbitrarily fix $\x\in\DmOl$ and particle volume decomposition $\PvDmSet=\{\PvDm{\i}\mid\i=1,2,\dots,\N\}$, and split $\errorfuncC{\intA}$ into
	\begin{align}
		\errorfuncC{\intA}(\x) & = \prffuncA{7}(\x) + \prffuncA{8}(\x) + \prffuncA{9}(\x) + \prffuncA{10}(\x)
	\end{align}	
	with
	\begin{align}
		\prffuncA{7}(\x) &\deq\errorfuncC{\intA}(\x)- \sum_{\i=1}^{\N}\sum_{\j=1}^{\N} \Meas{\VoroCell{\j}\cap\PvDm{\i}} (\x-\Pt{\i})^{\intA}|\wh(|\x-\Pt{\j}|)|, 
		\\
		\prffuncA{8}(\x) &\deq \sum_{\i=1}^{\N}\sum_{\j=1}^{\N} (\x-\Pt{\i})^{\intA} \int_{\VoroCell{\j}\cap\PvDm{\i}}\BrM{|\wh(|\x-\Pt{\j}|)|-|\wh(|\x-\y|)|}\dy,
		\\
		\prffuncA{9}(\x) &\deq \sum_{\i=1}^{\N}\sum_{\j=1}^\N \int_{\VoroCell{\j}\cap\PvDm{\i}}\{(\x-\Pt{\i})^{\intA} - (\x-\y)^{\intA}\}  |\wh(|\x-\y|)| \dy, 
		\\
		\prffuncA{10}(\x) &\deq \int_{\dRd}(\x-\y)^{\intA}|\wh(|\x-\y|)|\dy. 
	\end{align}	
	Then, we estimate $\prffuncA{7}$, $\prffuncA{8}$, $\prffuncA{9}$, and $\prffuncA{10}$. 
	
	From \eqref{sum_|VoroCell_PvDm|}, we can rewrite $\prffuncA{7}$ as
	\begin{equation}
		\prffuncA{7}(\x)=\sum_{\i=1}^{\N}\sum_{\j=1}^{\N} \Meas{\VoroCell{\j}\cap\PvDm{\i}}(\x-\Pt{\i})^{\intA} \{|\wh(|\x-\Pt{\i}|)|-|\wh(|\x-\Pt{\j}|)|\}. 
	\end{equation}
	For all $\y\in\DmH$, we have
	\begin{equation}
		|(\x-\y)^{\intA}|\leq |\x-\y|^{\intA}\leq \diam(\DmH)^{\intA}. 
		\label{prf:estimate:E7:01}
	\end{equation}
	From \eqref{prf:estimate:E1:02} and \eqref{prf:estimate:E7:01}, we obtain
	\begin{align}
		\BrA{\prffuncA{7}(\x)} 
		&\leq\sum_{\i=1}^{\N}\sum_{\j=1}^{\N}\Meas{\VoroCell{\j}\cap\PvDm{\i}}(\x-\Pt{\i})^{\intA}\big||\wh(|\x-\Pt{\i}|)|-|\wh(|\x-\Pt{\j}|)|\big|
		\nn\\
		&\leq\gconst\sum_{\i=1}^{\N}\sum_{\j=1}^{\N}\Meas{\VoroCell{\j}\cap\PvDm{\i}}\big||\wh(|\x-\Pt{\i}|)|-|\wh(|\x-\Pt{\j}|)|\big|
		\nn\\
		&\leq\gconst\sum_{\i=1}^{\N}\sum_{\j=1}^{\N}\Meas{\VoroCell{\j}\cap\PvDm{\i}}\BrA{\wh(|\x-\Pt{\i}|)-\wh(|\x-\Pt{\j}|)}
		\nn\\
		&\leq\gconst\BrS{1+\frac{\IndPtSet}{\h}}^{\Dim} \frac{\IndPvSetLocal{\PvDmSet}}{\h}. 
	\end{align}
	From \eqref{prf:estimate:E2:01} and \eqref{prf:estimate:E7:01}, we obtain
	\begin{align}
		|\prffuncA{8}(\x)| 
		&\leq \sum_{\i=1}^{\N}\sum_{\j=1}^{\N} (\x-\Pt{\i})^{\intA} \int_{\VoroCell{\j}\cap\PvDm{\i}}\big||\wh(|\x-\Pt{\j}|)|-|\wh(|\x-\y|)|\big|\dy
		\nn\\
		&\leq\gconst\sum_{\i=1}^{\N}\sum_{\j=1}^{\N}\int_{\VoroCell{\j}\cap\PvDm{\i}}\big||\wh(|\x-\Pt{\j}|)|-|\wh(|\x-\y|)|\big|\dy
		\nn\\
		&\leq\gconst\sum_{\i=1}^{\N}\sum_{\j=1}^{\N}\int_{\VoroCell{\j}\cap\PvDm{\i}}\BrA{\wh(|\x-\Pt{\j}|)-\wh(|\x-\y|)}\dy
		\nn\\
		& \leq \gconst \BrS{1+2\frac{\IndPtSet}{\h}}^{\Dim}\frac{\IndPtSet}{\h}.  
	\end{align}
	For all $\Pt\i\in\PtSet$ and $\y\in\DmH$, we have
	\begin{align}
		\BrA{(\x-\Pt{\i})^{\intA} - (\x-\y)^{\intA}}
		&=\BrA{\{(\x-\Pt{\i})-(\x-\y)\} \sum_{k=1}^{\intA} (\x-\Pt{\i})^{k-1}(\x-\y)^{\intA-k}}\nn\\
		&\leq\intA\,\diam(\DmH)^\intA\BrA{\y-\Pt{\i}}. 
		\label{prf:estimate:E9:01}
	\end{align}
	From \eqref{prf:estimate:E3:04}, \eqref{prf:estimate:E9:01}, and $\h<\H$, we obtain
	\begin{align}
		|\prffuncA{9}(\x)| 
		&\leq\sum_{\i=1}^{\N}\sum_{\j=1}^\N \int_{\VoroCell{\j}\cap\PvDm{\i}}\BrA{(\x-\Pt{\i})^{\intA} - (\x-\y)^{\intA}}  |\wh(|\x-\y|)| \dy
		\nn\\
		&\leq\gconst\sum_{\i=1}^{\N}\sum_{\j=1}^{\N} \int_{\VoroCell{\j}\cap\PvDm{\i}} |\y-\Pt{\i}| |\wh(|\x-\y|)| dy 
		\nn\\
		&\leq\gconst\BrS{1+\dfrac{\IndPtSet}{\h}}^{\Dim}\BrS{\IndPtSet+\IndPvSetLocal{\PvDmSet}}.
	\end{align}
	From \eqref{def:wh}, we obtain
	\begin{align}
		|\prffuncA{10}(\x)|
		&=\int_{\dRd}|\y|^{\intA}|\wh(|\y|)|\dy=\h^{\intA}\int_{\dRd}|\y|^{\intA}|\w(|\y|)|\dy. 
	\end{align}
	From the estimates of $\prffuncA{7}$, $\prffuncA{8}$, $\prffuncA{9}$, and $\prffuncA{10}$, we obtain
	\begin{equation}
		\NormCz{\errorfuncC{\intA}}{\DmOl} \leq \gconst\BrM{\BrS{1+2\frac{\IndPtSet}{\h}}^\Dim \dfrac{\IndPtSet+\IndPvSetLocal{\PvDmSet}}{\h}+\h^{\intA}}. 
	\end{equation}
	Because $\PvDmSet$ is arbitrary, we establish \eqref{eq:lem:Km}. 
\end{proof}

Using the lemmas defined above, we now prove Theorem \ref{thm:terror_cont_norm}. 
\begin{proof}[Theorem \ref{thm:terror_cont_norm}]
	By Lemmas \ref{lem:terror_cont_norm}, \ref{lem:Ja}, and \ref{lem:Km}, we have for all $\f\in\FsC{\WeightOrder+1}{\DmHOl}$ 
	\begin{equation}
		\NormCz{\f-\IntpApp{} \f}{\DmOl} \leq \gconst\BrM{ \BrS{1+2\frac{\IndPtSet}{\h}}^{\Dim}\frac{\IndPtSet+\IndPvSet}{\h}+\hN^{\WeightOrder+1}}\NormC{\f}{\WeightOrder+1}{\DmHOl}.
		\label{prf:ineq:terror_cont_norm_intp:01}
	\end{equation}
	Moreover, by Lemmas \ref{lem:terror_cont_norm}, \ref{lem:Jam}, and \ref{lem:Km}, when $\w$ satisfies Hypothesis \ref{hypo:ref_weight:smooth} with $\WeightOrderSmooth=0$, we have for all $\f\in\FsC{\WeightOrder+2}{\DmHOl}$ 
	\begin{equation}
		\NormCz{\Grad\f-\GradApp{}\f}{\DmOl}
		\leq \gconst\BrM{ \BrS{1+2\frac{\IndPtSet}{\h}}^{\Dim}\frac{\IndPtSet+\IndPvSet}{\h}+\hN^{\WeightOrder+1}}\NormC{\f}{\WeightOrder+2}{\DmHOl},
		\label{prf:ineq:terror_cont_norm_grad:01}
	\end{equation}
	and when $\w$ satisfies Hypothesis \ref{hypo:ref_weight:smooth} with $\WeightOrderSmooth=1$ for all $\f\in\FsC{\WeightOrder+3}{\DmHOl}$, 
	\begin{equation}
		\NormCz{\Lap\f-\LapApp{}\f}{\DmOl}
		\leq \gconst\BrM{ \BrS{1+2\frac{\IndPtSet}{\h}}^{\Dim}\frac{\IndPtSet+\IndPvSet}{\h^2}+\hN^{\WeightOrder+1}}\NormC{\f}{\WeightOrder+3}{\DmHOl}.
		\label{prf:ineq:terror_cont_norm_lap:01}
	\end{equation}
	Because the family $\{(\PtSet, \PvSet, \hN)\}_{\N\ra\infty}$ is regular, by applying \eqref{def:regular:cond} to \eqref{prf:ineq:terror_cont_norm_intp:01}, \eqref{prf:ineq:terror_cont_norm_grad:01}, and \eqref{prf:ineq:terror_cont_norm_lap:01}, we obtain \eqref{ineq:terror_cont_norm_intp}, \eqref{ineq:terror_cont_norm_grad}, and \eqref{ineq:terror_cont_norm_lap}, respectively.
	We now conclude the proof of Theorem \ref{thm:terror_cont_norm}. 
\end{proof}

\section{Numerical results}
\label{sec:Numerical_results}
\begin{figure}[tb]
	\includegraphics[width=50mm,bb=0 0 126.5mm 122.4mm]{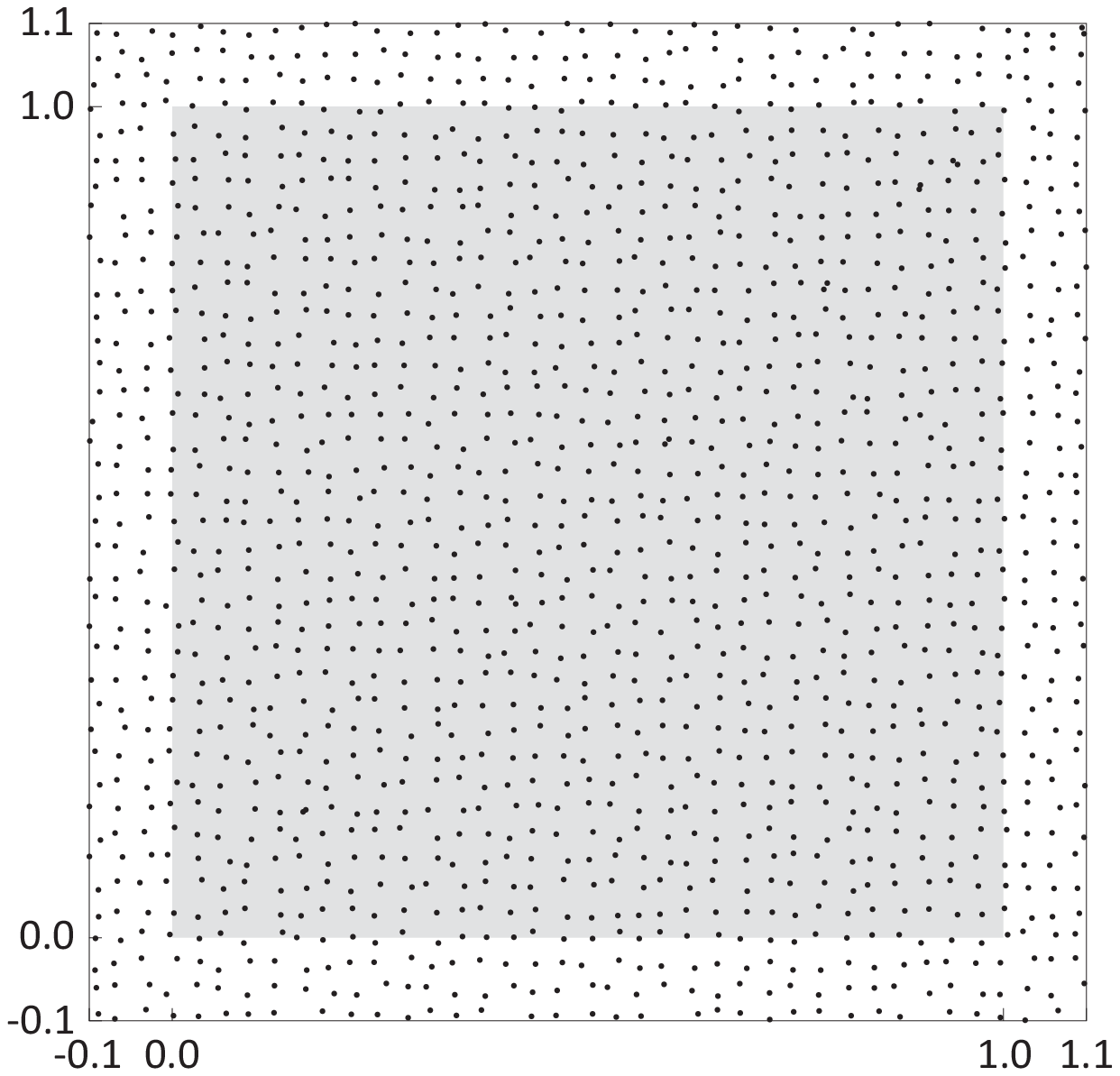}
	\caption{Particle distribution $\PtSet$ with $\Dx=2^{-5}$ $(\N=1,521)$. The gray area represents $\Dm$. }
	\label{fig_numerical_result_particle_distribution}
\end{figure}

Set $\Dm=(0,1)^2$ and $\H=0.1$. 
Then, $\DmH=(-0.1,1.1)^2$. 
We now compute the truncation errors of $\f:\DmH\ra\dR$, which are defined as $\f(x,y)=\sin(2\pi(x+y))$. 
Particle distribution $\PtSet$ is set as
\begin{equation}
	\PtSet = \BrM{\BrS{(i+\eta_{ij}^{(1)})\Delta x, (j+\eta_{ij}^{(2)}) \Delta x}\in\DmH ;~i,j\in\dZ},
\end{equation}
where $\Dx$ is taken by $2^{-5}, 2^{-6}, \dots, 2^{-12}$ and $\Noise{\i}{\j}{k}\,(i,j\in\dZ, k=1,2)$ are random numbers satisfying $|\Noise{\i}{\j}{k}|< 1/4$. 
Particle distribution $\PtSet$ with $\Dx=2^{-5}$ is shown in Figure \ref{fig_numerical_result_particle_distribution}.
Particle volume set $\PvSet$ is defined as
\begin{equation}
	\PvSet=\BrM{\Pv{\i} = \dfrac{\Meas{\DmH}}{\N}\midd \i=1,2,\dots,\N}.  
	\label{pvolume_exp_02} 
\end{equation}
For $m=1,3,5$, the influence radius $\hN$ is set as 
\begin{align}
	\hN = 2.6\times2^{5/\RegOrder-5}\Dx^{1/\RegOrder}. 
\end{align}
Note that if $\Dx=2^{-5}$, then $\h=2.6\times2^{-5}$ for all $m$.  
Using the discrete parameters above, 
the covering radius $\IndPtSet$ satisfies $\IndPtSet \leq \sqrt{2}(1+1/4)\Dx/2$. 
Moreover, 
the Voronoi deviation $\IndPvSet$ satisfies $\IndPvSet\leq64(1+\sqrt{2})\Dx/\pi$. 
Therefore, the family $\{(\PtSet,\PvSet,\hN)\}$ is regular with order $m$. 

For the interpolant, we consider the following three cases of reference weight functions: 
\begin{align*}
	({\rm \Intp1})\quad \w(r) &\deq\dfrac{3}{\pi} 
	\begin{cases}
		1-r,& 0 \leq r < 1,
		\\
		0,\hspace{16ex}& 1 \leq r, 
	\end{cases}
	\\
	({\rm \Intp2})\quad \w(r) &\deq \dfrac{40}{7\pi}
	\begin{cases}
		1-6r^2+6r^3, \quad & \ds 0 \leq r < \frac{1}{2},
		\vspace*{0.5ex}\\
		2(1-r)^3, \quad & \ds \frac{1}{2} \leq r < 1,
		\vspace*{0.5ex}\\
		0, \hspace{13.5ex} & 1 \leq r, 
	\end{cases}
	\\
	({\rm \Intp3})\quad\w(r) &\deq\dfrac{5}{\pi} 
	\begin{cases}
		(1-r)(2-3r), & 0 \leq r < 1,
		\vspace*{1ex}\\
		0,\hspace{16ex}& 1 \leq r. 
	\end{cases}
\end{align*}
$({\rm \Intp1})$ is the lowest order polynomial function belonging to $\FsWeightFunc$. 
$({\rm \Intp2})$ is the cubic B-spline commonly used in SPH and belonging to $\FsWeightFunc$. 
$({\rm \Intp3})$ is the lowest order polynomial function belonging to $\FsWeightFunc$ that satisfies Hypothesis \ref{hypo:ref_weight:order} with $\WeightOrder=3$.  

For the approximate gradient operator, we consider the following three cases of reference weight functions: 
\begin{align*}
	({\rm \Grad1})\quad \w(r) &\deq\dfrac{6}{\pi}
	\begin{cases}
		r(1-r), \quad & 0 \leq r < 1,
		\\
		0,\hspace{19ex}& 1 \leq r, 
	\end{cases}
	\\
	({\rm \Grad2})\quad \w(r) &\deq \dfrac{40}{7\pi}
	\begin{cases}
		6r^2-9r^3, \quad & \ds 0 \leq r < \frac{1}{2},
		\vspace*{0.5ex}\\
		3r(1-r)^2, \quad & \ds \frac{1}{2} \leq r < 1,
		\vspace*{0.5ex}\\
		0,\hspace{17.5ex} & 1 \leq r, 
	\end{cases}
	\\
	({\rm \Grad3})\quad \w(r) &\deq\dfrac{15}{2\pi} 
	\begin{cases}
		r(1-r)(5-7r), & 0 \leq r < 1,
		\\
		0,\hspace{18ex}& 1 \leq r. 
	\end{cases}
\end{align*}
$({\rm \Grad1})$ is the lowest order polynomial function belonging to $\FsWeightFunc$ and Hypothesis \ref{hypo:ref_weight:smooth} with $\WeightOrderSmooth=0$. 
$({\rm \Grad2})$ is chosen so that approximate gradient operator \eqref{gh_def} with $({\rm \Grad2})$ coincides with that in SPH with the cubic B-spline (see Appendix \ref{sec:appendix_conventional_particle_methods}). 
$({\rm \Grad3})$ is the lowest order polynomial function belonging to $\FsWeightFunc$, Hypothesis \ref{hypo:ref_weight:order} with $\WeightOrder=3$, and Hypothesis \ref{hypo:ref_weight:smooth} with $\WeightOrderSmooth=0$.  

For the approximate Laplace operator, we consider the following three cases of reference weight functions: 
\begin{align*}
	({\rm \Lap1})\quad \w(r) &\deq\dfrac{10}{\pi}
	\begin{cases}
		r^2(1-r), \quad & 0 \leq r < 1,
		\\
		0,\hspace{21ex}& 1 \leq r, 
	\end{cases}
	\\
	({\rm \Lap2})\quad \w(r) &\deq \dfrac{40}{7\pi}
	\begin{cases}
		6r^2-9r^3, \quad & \ds 0 \leq r < \frac{1}{2},
		\vspace*{0.5ex}\\
		3r(1-r)^2, \quad & \ds \frac{1}{2} \leq r < 1,
		\vspace*{0.5ex}\\
		0,\hspace{20.5ex} & 1 \leq r, 
	\end{cases}
	\\
	({\rm \Lap3})\quad \w(r) &\deq\dfrac{30}{\pi} 
	\begin{cases}
		r^2(1-r)(3-4r), & 0 \leq r < 1,
		\\
		0,\hspace{21ex}& 1 \leq r. 
	\end{cases}
\end{align*}
$({\rm \Lap1})$ is the lowest order polynomial function belonging to $\FsWeightFunc$ and Hypothesis \ref{hypo:ref_weight:smooth} with $\WeightOrderSmooth=1$. 
$({\rm \Lap2})$ is chosen so that approximate Laplace operator \eqref{lh_def} with $({\rm \Lap2})$ coincides with that in SPH with the cubic B-spline (see Appendix \ref{sec:appendix_conventional_particle_methods}). 
$({\rm \Lap3})$ is the lowest order polynomial function belonging to $\FsWeightFunc$, Hypothesis \ref{hypo:ref_weight:order} with $\WeightOrder=3$, and Hypothesis \ref{hypo:ref_weight:smooth} with $\WeightOrderSmooth=1$.  

\begin{figure}[h!]
	\includegraphics[width=100mm,bb=0 0 186.0mm 286.0mm]{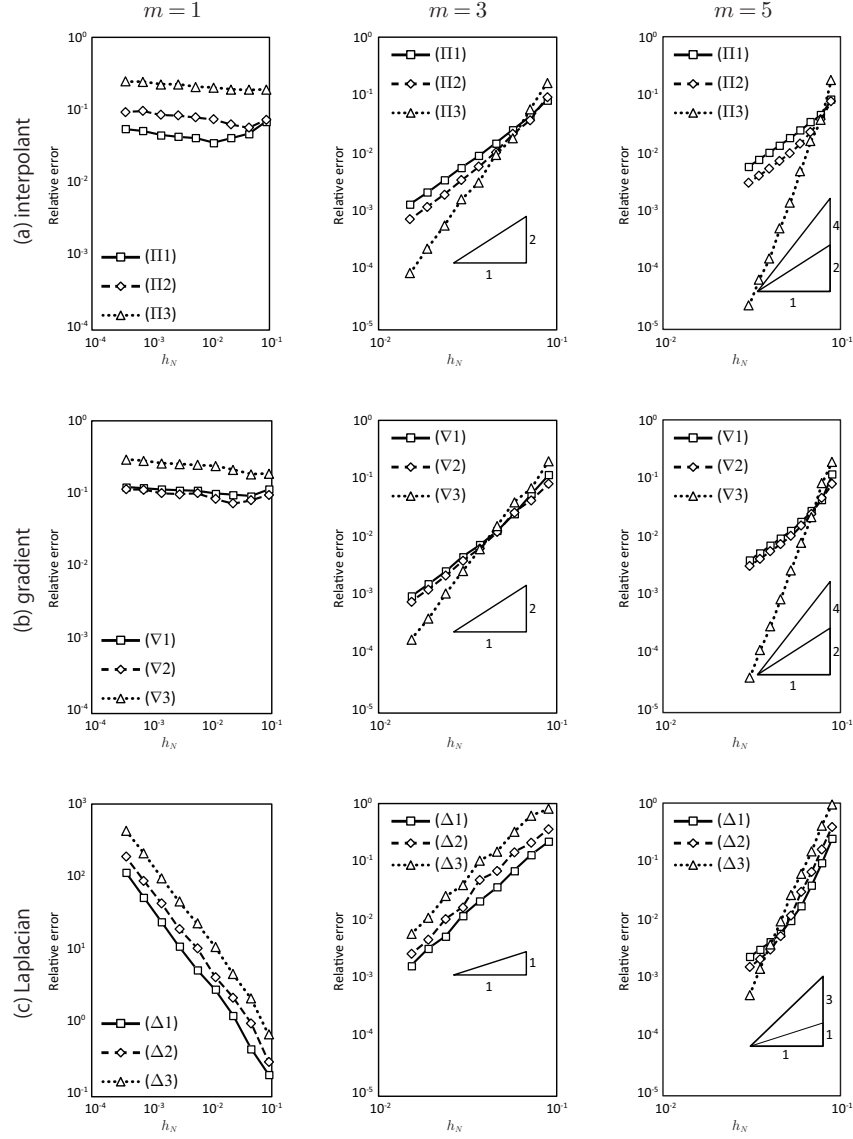}
	\caption{Graphs of the relative errors of (a) the interpolant, (b) approximate gradient operator, and (c) approximate Laplace operator versus the influence radius with regular orders $m=1, 3, 5$.}
	\label{fig:numerical_terror}
\end{figure}

The above settings were used in the computation of the following relative errors
\begin{equation}
	\dfrac{\NormDiscL{\f-\IntpApp{} \f}{\infty}{\Dm}}{\NormCzNd{\f}{\DmOl}},\quad
	\dfrac{\NormDiscL{\Grad\f-\GradApp{} \f}{\infty}{\Dm}}{\NormCzNd{\Grad\f}{\DmOl}},\quad
	\dfrac{\NormDiscL{\Lap\f-\LapApp{} \f}{\infty}{\Dm}}{\NormCzNd{\Lap\f}{\DmOl}}. 
\end{equation}
Here, the discrete norm $\NormDiscL{\cdot}{\infty}{\Dm}$ is defined as
\begin{align}
	\NormDiscL{\f}{\infty}{\Dm} \deq \max_{\i\in\IndexSet{\Dm}} |\f(\Pt{\i})|. 
\end{align}
Figure \ref{fig:numerical_terror} shows graphs of the relative errors of (a) interpolant $\IntpApp{}$, (b) approximate gradient operator $\GradApp{}$, and (c) approximate Laplace operator $\LapApp{}$ versus the influence radius $\hN$ with regular orders $\RegOrder=1, 3, 5$. 
In Figure \ref{fig:numerical_terror}, the slopes of the triangles show the theoretical convergence rates obtained via Theorem \ref{thm:terror_cont_norm}. 
Table \ref{tab:numerical_theoretical_converge_rate} lists the numerical and theoretical convergence rates obtained from the cases of $\Delta x=2^{-11}$ and $2^{-12}$, where the theoretical convergence rates correspond to Theorem \ref{thm:terror_cont_norm}. 
In the case of $\RegOrder=1$, as the settings could not be applied to Theorem \ref{thm:terror_cont_norm}, only numerical results without convergence were obtained. 
In contrast, the settings in cases $\RegOrder=3$ and $5$ could be applied Theorem \ref{thm:terror_cont_norm}; thus, the numerical results with convergence were obtained. 
Moreover, the approximate operators with reference weight functions satisfying Hypothesis \ref{hypo:ref_weight:order} with $\WeightOrder=3$ became higher convergence orders in the cases where $\RegOrder=5$ as per Theorem \ref{thm:terror_cont_norm}. 
\begin{table}[htb]
	\renewcommand{\arraystretch}{1.2}
	\caption{Numerical and theoretical convergence rates of (a) the interpolant, (b) approximate gradient operator, and (c) approximate Laplace operator with regular orders $m=1, 3, 5$. The numerical convergence rates were obtained for the cases of $\Delta x=2^{-11}$ and $2^{-12}$. }
	\label{tab:numerical_theoretical_converge_rate}
	\begin{tabular}{lrrrrrr} 
		\hline
		&\multicolumn{2}{c}{$m=1$}&\multicolumn{2}{c}{$m=3$}&\multicolumn{2}{c}{$m=5$}\\
		& Numer. & Theor. & Numer. & Theor. & Numer. & Theor.
		\\\hline
		\noalign{\hrule height 1pt}
		\\[-2.5ex]
		(a) Interpolant
		\\[0.5ex]
		($\rm \Intp1$) & -0.10 & N/A &2.02 & 2 & 2.00 & 2\\
		($\rm \Intp2$) & 0.05 & N/A & 2.13 & 2 & 2.01 & 2\\
		($\rm \Intp3$) &  0.00 & N/A & 4.20 & 2 & 7.41 & 4
		\\[1.5ex]
		(b) Gradient
		\\[0.5ex]
		($\rm \Grad1$) & -0.05 & N/A & 2.11 & 2 & 2.02 & 2\\
		($\rm \Grad2$) & -0.02 & N/A & 2.08 & 2 & 2.03 & 2\\
		($\rm \Grad3$) &  -0.06 & N/A & 3.56 & 2 & 7.69 & 4
		\\[1.5ex]
		(c) Laplacian
		\\[0.5ex]
		($\rm \Lap1$) & -1.14 & N/A &2.91 & 1 & 2.05 & 2\\
		($\rm \Lap2$) & -1.09 & N/A & 2.39 & 1 & 2.23 & 2\\
		($\rm \Lap3$) & -1.02 & N/A & 2.76 & 1 & 7.50 & 3\\
		\hline
	\end{tabular}
\end{table}

\section{Conclusions}
\label{sec:conclusions}
We analyzed truncation errors in a generalized particle method, which is a wider class of particle methods that includes commonly used methods such as SPH and MPS. 
In our analysis, we introduced two indicators: 
the first was the covering radius, which represents the maximum radius of the Voronoi region associated with the particle distribution, while the second was the Voronoi deviation, which indicates the deviation between particle volumes and Voronoi volumes. 
With the covering radius and Voronoi deviation, we introduced a regularity of a family of discrete parameters, which includes the particle distribution, particle volume set, and influence radius associated with the number of particles. 
Moreover, we introduced two hypotheses of reference weight functions. 
With the regularity and hypotheses of reference weight functions, we established truncation error estimates for the continuous norm. 
The convergence rates are dependent on the regular order and order of the reference weight functions appearing in a hypothesis. 
Moreover, as it was possible to validate the conditions by calculation, we showed the numerical convergence orders were in good agreement with the theoretical ones. 

In a forthcoming paper, we plan to establish error estimates of the generalized particle method for the Poisson and heat equations.

\section*{Acknowledgments}
We would like to thank Drs. Daisuke Tagami and Hayato Waki for their helpful comments over the course of this work.
This research was supported by JSPS KAKENHI Grant Number 17K17585 and the JSPS A3 Foresight Program. 

\appendix
\section{Description of conventional particle methods by the generalized particle method}
\label{sec:appendix_conventional_particle_methods}
This appendix provides a description of conventional particle methods, such as the smoothed particle hydrodynamics (SPH) \cite{liu2010smoothed,price2012smoothed} and the moving particle semi-implicit (MPS) \cite{koshizuka1996moving}, in the context of the generalized particle method. 
In SPH, by using the reference weight function $\ws\in\FsWeightFunc$ and parameters $m_\i, \rho_\i\in\dRp\,(\i=1,2,\dots,\N)$, for $\f\in\FsCz{\DmHOl}$, the approximate operators are defined as
\begin{align}
	\IntpAppSPH \f(\x) &\deq \sum_{\i=1}^\N\dfrac{m_\i}{\rho_\i}\f(\Pt{\i})\whs(|\x-\Pt{\i}|),\qquad\x\in\DmH,
	\label{ih_SPH_def}\\
	\GradAppSPH \f(\x) &\deq \sum_{\i=1}^\N\dfrac{m_\i}{\rho_\i} \BrM{\f(\x)-\f(\Pt{\i})}\nabla\whs(|\x-\Pt{\i}|),\qquad\x\in\DmH,
	\label{gh_SPH_def}\\
	\LapAppSPH \f(\x) &\deq 2 \sum_{\i\in\IndexSetXHAst{\x}{\h}}\dfrac{m_\i}{\rho_\i} \frac{\f(\x)-\f(\Pt{\i})}{|\x-\Pt{\i}|}\frac{\x-\Pt{\i}}{|\x-\Pt{\i}|}\cdot\nabla\whs(|\x-\Pt{\i}|),\qquad\x\in\DmH. 
	\label{lh_SPH_def}
\end{align}
By setting $\w=\ws$ and $\PvSet=\SetNd{\Pv{\i}=m_\i/\rho_\i}{\i=1,2,\dots,\N}$, the generalized interpolant \eqref{ih_def} coincides with \eqref{ih_SPH_def}. 
Moreover, because 
\begin{align}
	- \int_{\dRd} \dfrac{\x}{\Dim}\cdot \nabla \ws(|\x|) \dx = \int_{\dRd} \ws(|\x|) \dx = 1,  
\end{align}
by setting
\begin{equation}
	\w(r)=-\Dim^{-1} r \frac{d}{dr} \ws(r),
	\label{weight_general_SPH}
\end{equation}
and $\PvSet=\SetNd{\Pv{\i}=m_\i/\rho_\i}{\i=1,2,\dots,\N}$, \eqref{gh_def} and \eqref{lh_def} coincide with \eqref{gh_SPH_def} and \eqref{lh_SPH_def}, respectively. 

From Theorem \ref{thm:terror_cont_norm}, we obtain the following corollary that is a truncation error estimate of approximate operators \eqref{gh_SPH_def} and \eqref{lh_SPH_def}. 
\begin{corollary}
	\label{rem:terror_cont_norm_SPH}
	Suppose that parameters $\rho_\i, m_\i$ satisfy
	\begin{equation}
		\sum_{\i=1}^\N\dfrac{m_\i}{\rho_\i} = \Meas{\DmH}, 
	\end{equation}
	and that $\{(\PtSet, \PvSet, \hN)\}_{\N\ra\infty}$ is regular with order $\RegOrder$, where $\PvSet=\SetNd{\rho_\i/m_\i}{\i=1,2,\dots,\N}$. 
	Moreover, suppose that $\ws$ satisfies the following conditions{\rm ;} 
	\begin{equation}
		\ws\in\FsC{2}{\dRpz}, \quad
		\dfrac{d}{dr} \ws(r) < 0~(0< r< 1),\quad
		\lim_{s\downarrow 0}\BrA{\dfrac{1}{s} \dfrac{d}{dr} \ws(s)} < \infty. 
		\label{ws:cond}
	\end{equation}
	Then, there exists a positive constant $\thmconst$ independent of $\N$ such that
	\begin{alignat}{3}
		\NormCz{\f-\IntpAppSPH\f}{\DmOl} & \leq \thmconst\, \h^{\min\{2,\RegOrder-1\}} \NormC{\f}{2}{\DmHOl},\qquad \f\in \FsC{2}{\DmHOl}, 
		\label{eq:terror_cont_norm_intp_SPH}
		\\
		\NormCz{\Grad\f-\GradAppSPH\f}{\DmOl} & \leq \thmconst\, \h^{\min\{2,\RegOrder-1\}} \NormC{\f}{3}{\DmHOl},\qquad \f\in \FsC{3}{\DmHOl}, 
		\label{eq:terror_cont_norm_grad_SPH}
		\\
		\NormCz{\Lap\f-\LapAppSPH\f}{\DmOl} & \leq \thmconst\, \h^{\min\{2,\RegOrder-2\}} \BrN{\f}_{\FsC{4}{\DmHOl}},\qquad \f\in \FsC{4}{\DmHOl}. 
		\label{eq:terror_cont_norm_lap_SPH}
	\end{alignat}
\end{corollary}

\begin{remark}
	Note that representative reference weight functions employed in SPH, such as the cubic B-spline, quintic B-spline, and Wendland function {\rm(}5-order positive definite function{\rm)} {\rm \cite{dehnen2012improving,liu2010smoothed}}, satisfy \eqref{ws:cond}. 
\end{remark}

In MPS \rmn{\cite{khayyer2011enhancement}}, by using reference weight function $\wm\in\FsWeightFunc$ 
and parameters $\ParamMPSNZero, \ParamMPSLambda\in\dRp$ for $\f\in\FsC{1}{\DmHOl}$, approximate differential operators can be defined as
\begin{align}
	\GradAppMPS \f(\x) &\deq \frac{\Dim}{\ParamMPSNZero} \sum_{\i\neq\j} \frac{\f(\Pt{\i})-\f(\x)}{|\x-\Pt{\i}|}\frac{\Pt{\i}-\x}{|\x-\Pt{\i}|}\whm(|\x-\Pt{\j}|),\qquad\x\in\DmH,
	\label{gh_MPS_def}\\
	\LapAppMPS \f(\x) &\deq \frac{2\Dim}{\ParamMPSNZero\ParamMPSLambda} \sum_{\i\neq\j} \BrM{\f(\Pt{\i})-\f(\x)}\whm(|\x-\Pt{\j}|),\qquad\x\in\DmH. 
	\label{lh_MPS_def}
\end{align}
Note that an interpolant is not defined in MPS. 
By setting $\w=\wm$ and $\PvSet=\SetNd{\Pv{\i}=\ParamMPSNZero^{-1}}{\i=1,2,\dots,\N}$, the approximate gradient operator \eqref{gh_def} coincides with \eqref{gh_MPS_def}. 
Moreover, by setting $\w(r) = \ParamMPSLambda^{-1}r^2\wm(r)$ and $\PvSet=\SetNd{\Pv{\i}=\ParamMPSNZero^{-1}}{\i=1,2,\dots,\N}$, 
approximate Laplace operator \eqref{lh_def} coincides with \eqref{lh_MPS_def}. 

\begin{corollary}
	\label{cor:terror_MPS}
	Suppose that
	\begin{equation}
		\ParamMPSNZero=\dfrac{\N}{\Meas{\DmH}},
		\qquad
		\ParamMPSLambda=\int_{\dRd}|\x|^2\wm(|\x|)\dx,  
		\qquad 
		\wm\in\FsWeightFunc. 
	\end{equation}
	Moreover, suppose that $\{(\PtSet, \PvSet, \hN)\}_{\N\ra\infty}$ is regular with order $\RegOrder$, where $\PvSet=\SetNd{\Pv{\i}=\ParamMPSNZero^{-1}}{\i=1,2,\dots,\N}$. 
	Then, there exists a positive constant $\thmconst$ independent of $\N$ such that
	\begin{alignat}{3}
		\NormCz{\Lap\f-\LapAppMPS\f}{\DmOl} & \leq \thmconst\, \h^{\min\{2,\RegOrder-2\}} \BrN{\f}_{\FsC{4}{\DmHOl}},\qquad \f\in \FsC{4}{\DmHOl}. 
		\label{eq:terror_cont_norm_lap_MPS}
	\end{alignat}
	Furthermore, when $\wm$ satisfies Hypothesis \ref{hypo:ref_weight:smooth} with $\WeightOrderSmooth=0$, 
	\begin{alignat}{3}
		\NormCz{\Grad\f-\GradAppMPS\f}{\DmOl} & \leq \thmconst\, \h^{\min\{2,\RegOrder-1\}} \NormC{\f}{3}{\DmHOl},\qquad \f\in \FsC{3}{\DmHOl}. 
		\label{eq:terror_cont_norm_grad_MPS}
	\end{alignat}
\end{corollary}

\begin{remark}
	Note that the reference weight function, which is commonly used in the MPS and defined as
	\begin{align}
		\wm(r) &\deq 
		\begin{cases}
			\dfrac{1}{r}-1, \qquad & \ds 0 \leq r < 1,
			\smallskip\\
			0, \qquad & 1 \leq r, 
		\end{cases}
		\label{w_MPS}
	\end{align}
	does not satisfy $\wm\in\FsWeightFunc$. 
	In contrast, the continuous reference weight function as introduced in \cite{shakibaeinia2010weakly} satisfies $\wm\in\FsWeightFunc$. 
	However, as far as we know, no reference weight functions that also satisfy Hypothesis \ref{hypo:ref_weight:smooth} with $\WeightOrderSmooth=0$ are proposed in MPS.
\end{remark}

\section{Computational procedure of the indicators}
\label{sec:appendix_compt_procedure_indicator}
This appendix introduces the procedures for computing the indicators introduced in this paper, namely, the covering radius \eqref{def:covering_radius} and Voronoi deviation \eqref{def:Voronoi_deviation}. 

The covering radius $\IndPtSet$ can be computed as follows. 
As per the methods used to construct Voronoi decompositions, such as the increment method \citerm{boissonnat1998algorithmic}, we first draw the boundaries of the Voronoi region in $\DmH$. 
Next, for each particle, we compute the maximum distance from particle $\Pt{\i}$ to the boundary of its Voronoi region $\VoroCell{\i}$ (i.e., $\max_{\y\in\Ol{\VoroCell{\i}}}|\Pt{\i}-\y|$). 
Finally, we obtain the covering radius $\IndPtSet$ by computing the maximum of these distances. 

Next, we consider the Voronoi deviation $\IndPvSet$. 
Let $\zeta\in\dR^{3\N}$ be 
\begin{align}
	\zeta &\deq(\Meas{\VoroCell{1}},\Meas{\VoroCell{2}},\dots,\Meas{\VoroCell{\N}},\Pv{1},\Pv{2},\dots,\Pv{\N},0,0,\dots,0)^T. 
\end{align}
Using parameters $q, s_i, a_{ij}\in\dRp$ $(i,j=1,2,\dots,\N)$, we set $z\in\dR^{N^2+N+1}$ as
\begin{equation}
	z\deq(a_{11},a_{12},\dots,a_{NN},s_1,s_2,\dots,s_N,q)^T. 
\end{equation}
Moreover, we set $M\in\dR^{3\N\times(\N^2+\N+1)}$ so that equation $Mz=\zeta$ represents 
\begin{align}
	\sum_{\j=1}^{\N} a_{ij} = \Meas{\VoroCell{\i}}, \qquad \sum_{\j=1}^{\N} a_{ji} = \Pv{\i}, \qquad i=1,2,\dots,N
	\label{cond_aij_01}
\end{align}
and
\begin{align}
	q=s_i+\sum_{\j=1}^{\N} \frac{a_{ij}+a_{ji}}{\Meas{\VoroCell{\i}}}|\Pt{\i}-\Pt{\j}|,\qquad i=1,2,\dots,N. 
	\label{cond_q_si}
\end{align}
Then, by considering $a_{ij}$ to be $\Meas{\VoroCell{\i}\cap\PvDm{\j}}$, we find that the minimum value of $q$ with condition $Mz=\zeta$ coincides with the Voronoi deviation $\IndPvSet$. 
We therefore consider the linear problem: 
\begin{align}
	\mbox{Minimize} \quad b^Tz\quad\mbox{subject to} \quad Mz=\zeta,~z\geq 0. 
	\label{optimize}
\end{align}
Here, $b\deq(0,0,\dots,0,1)^T\in\dR^{\N^2+\N+1}$. 
The solution $b^Tz$ of \eqref{optimize} is equivalent to the Voronoi deviation $\IndPvSet$. 
Because $Mz=\zeta$ is unique for $(\PtSet, \PvSet, \hN)$, the linear problem is computable via numerical methods for linear programming problems, such as the simplex method \citerm{dantzig1998linear}. 

\section{Construction of reference weight functions}
\label{sec:appendix_const_ref_weight_func}
For all $\WeightOrder\in\dN$ $(\WeightOrder\geq2)$, it is possible to construct a reference weight function satisfying Hypothesis \ref{hypo:ref_weight:order} with $\WeightOrder$ as the condition of Hypothesis \ref{hypo:ref_weight:order} can be rewritten to include a finite number of conditions
\begin{align}
	\int_{0}^1 r^{\Dim+2j-1} w(r) dr = 0,\qquad  j=1,2,\dots,\gausssymbol{\WeightOrder/2}. 
	\label{hypo:ref_weight:order_2}
\end{align}
Here, the Gauss symbol $\gausssymbol{a}$  denotes the largest integer that is less than or equal to $a$. 
For example, function $\w$ is set as the $p$th polynomial function: 
\begin{align}
	\w(r) \deq 
	\begin{cases}
		\ds \gamma_\Dim\BrS{1+\sum_{\ell=1}^{p} a_\ell r^\ell},\quad &0\leq r<1,
		\\
		\ds 0,\quad &r \geq 1.
	\end{cases}
\end{align}
Then, if coefficients $a_\ell$ satisfy the linear equations
\begin{align}
	\gamma_\Dim\BrS{1+\sum_{\ell=1}^{p} \dfrac{a_\ell}{\ell+\Dim}} &= 1,
	\\
	\sum_{\ell=1}^{p} a_\ell &= 0,
	\\
	\sum_{\ell=1}^{p} \ell a_\ell &= 0,
	\\
	1 + \sum_{\ell=1}^{p} \dfrac{\Dim+2j}{\Dim+\ell+2j}a_\ell &= 0,\qquad j=1,2,\dots,\gausssymbol{\WeightOrder/2}, 
\end{align}
then $\w$ satisfies $\w\in\FsWeightFunc$ and Hypothesis \ref{hypo:ref_weight:order} with $\WeightOrder$. 
Therefore, to construct reference functions with Hypothesis \ref{hypo:ref_weight:order} with $\WeightOrder$ represented by polynomial functions, the degree of the polynomial functions must be at least $\gausssymbol{\WeightOrder/2}+2$.  

Moreover, for all $\WeightOrderSmooth\in\dNz$, reference weight functions satisfying Hypothesis \ref{hypo:ref_weight:smooth} with $\WeightOrderSmooth$ can be constructed based on the following proposition. 

\begin{proposition}
	\label{prop_w_cond_01}
	Assume that reference the weight function $\w$ defined in $\dRpz$ satisfies $\w\in C^1(\dRpz)$ and is represented by a polynomial function in $[0,s]$ for $s\in(0,1]$. 
	Let $p_{0}$ be the minimum degree of $\w$ in $[0,s]$. 
	Then, if $p_{0}-\WeightOrderSmooth\geq1$, $\w$ satisfies Hypothesis \ref{hypo:ref_weight:smooth} with $\WeightOrderSmooth$. 
\end{proposition}

\begin{proof}
	From the assumption, $\w$ can be represented by 
	\begin{align}
		\w(r) =
		\sum_{\ell=p_0}^{p} a_\ell r^{\ell}, \qquad & 0\leq r< s,
	\end{align}
	where $p\in\dN$ and $a_\ell\in\dR$ $(\ell=p_0, p_0+1, \dots, p)$. 
	Set $\wn{\WeightOrderSmooth}$ as \eqref{def:wn}. 
	Since 
	\begin{align}
		\sup_{r\in(0,s)}\BrANd{\wn{\WeightOrderSmooth+1}(r)} &\leq\sum_{\ell=p_0}^{p} |a_\ell| s^{\ell-\WeightOrderSmooth-1}<\infty,
		\\
		\sup_{r\in(s,\infty)}\BrANd{\wn{\WeightOrderSmooth+1}(r)} &\leq\dfrac{\NormCz{\w}{[s,1]}}{s^{\WeightOrderSmooth+1}}<\infty
	\end{align}
	and 
	\begin{align}
		\sup_{r\in(0,s)}\BrA{\dfrac{d}{dr}\wn{\WeightOrderSmooth}(r)} &\leq\sum_{\ell=p_0}^{p} (\ell-\WeightOrderSmooth) |a_\ell| s^{\ell-\WeightOrderSmooth-1}<\infty,
		\\
		\sup_{r\in(s,\infty)}\BrA{\dfrac{d}{dr}\wn{\WeightOrderSmooth}(r)} &\leq\dfrac{\WeightOrderSmooth\NormC{\w}{1}{[s,1]}}{s^{\WeightOrderSmooth+1}}<\infty,
	\end{align}
	if $p_0-\WeightOrderSmooth\geq1$, we have $\w$ satisfies Hypothesis \ref{hypo:ref_weight:smooth} with $\WeightOrderSmooth$. 
\end{proof}
This proposition means that the regularity of the reference functions around zero is important.

\end{document}